\documentclass[11pt,a4paper]{article}
\usepackage{t1enc}
\usepackage[latin1]{inputenc}
\usepackage[english]{babel}
\usepackage{harvard}
\usepackage{amssymb,stmaryrd}

\newcommand{\C}{\mathbb C}
\newcommand{\R}{\mathbb R}

\newcommand{\E}{\mathbb E}

\newcommand{\beq}{\begin{equation}}
\newcommand{\eeq}{\end{equation}}
\newcommand{\beqarr}{\begin{eqnarray}}
\newcommand{\eeqarr}{\end{eqnarray}}
\newcommand{\beqa}{\begin{eqnarray*}}
\newcommand{\eeqa}{\end{eqnarray*}}
\newtheorem{theorem}{Theorem}

\newtheorem{lemma}{Lemma}
\newtheorem{result}{Result}

\begin{document}
\thispagestyle{empty}

\title{The problem of space  in the light of relativity: \\  the views of H. Weyl and E. Cartan }
\author{Erhard Scholz\footnote{University of Wuppertal, Department C,  Math., and IZWT,\newline  \hspace*{1.3em}  scholz@math.uni-wuppertal.de}}

\date{13. 07. 2014 }
\maketitle
%\vspace{2.5mm}

\begin{abstract}
Starting from a short review of the ``classical'' space problem in the sense of the 19th century (Helmholtz --  Lie -- Klein) it is discussed how the  challenges posed by special and general relativity to the classical analysis were taken up by Hermann Weyl and Elie Cartan. Both mathematicians reconsidered  the space problem from the point of view of  transformations operating in the infinitesimal neighbourhoods of a manifold (spacetime). In a short outlook we survey further developments in mathematics and physics of the second half of the 20th century, in which  core ideas of Weyl's and/or Cartan's analysis of the space problem were further investigated (mathematics) or incorporated into basic theories (physics).
\end{abstract}

\subsubsection*{1. Introduction}
Once Euclidean geometry was no longer considered the only possibility for describing physical space, the question arose how  the structure of  space, or of spacetime, could be characterized in terms of  symmetry structures on manifolds. The question was posed and discussed by Riemann (1854),   taken up anew  by Helmholtz at the end of the 1860s with ingenious conceptual insight but  mathematically quite vague, even from the standpoint of contemporary mathematics.  Helmholtz's argument was  refined   by Lie (1886ff.) and, independently, by Killing (1885ff.). It acquired a prominent place in the broader discourse on mathematics and reality at the turn to the 20th century. In these considerations   the possibility of ``freely moving'' rigid bodies as measuring devices was crucial. We shall speak of the {\em classical} analysis of the {\em problem of space} (PoS). 

At the beginning of the 20th century special and general relativity changed the conceptual scene and the interface with empirical reality drastically. Special relativity brought in the inseparable union of space and time and undermined the classical stipulation of  ``rigidity''. The general theory complicated things even further, because now the classical assumption of metrical transformations  in finite regions, or even globally, lost any ground.  The arising questions was under discussion among physicists, e.g.  \cite{Einstein:GeometrieundErfahrung},  philosophers   \cite{Schlick:Raum,Carnap:Raum,Reichenbach:1920,Reichenbach:1928}, and  mathematicians, in particular H. Weyl and E. Cartan. 
Weyl and Cartan had particular interest in the question of  how to reformulate the older criteria of the choice of space structure in the light of differential geometry and  the general theory of relativity ({\em modern} or {\em relativistic} PoS). They chose different roads and arrived at different conclusions, but there was an overlap and interchange between the two authors. 

Our paper  analyzes their  respective ways of posing the question and their path toward an answer (section 4 and 5).\footnote{For this aspect see also \cite{Scholz:Weyl_PoS,Scholz:EMS}.} Before we do so, we review  the classical problem of space and sketch how it was undermined by  the two relativity theories (sections 2 and 3). In a short survey (section 6) we follow more recent developments by which Weyl's and Cartan's infinitesimal geometric refinement of the analysis of the PoS fanned out in mathematics and physics of the second half of the 20th century.\footnote{A similar question regarding the spread of Klein's view in physics during the 20th century is being surveyed in \cite{Goenner:2013}.}  
The most important  technical  terms of more recent origin, not necessarily common to all readers interested in the history of the concept of space, are explained (loosely) in a glossary at the end of the article.

\subsubsection*{2. The classical problem of space}
Hermann Helmholtz (after 1883 {\em von}  Helmholtz)  based his analysis of the problem of space about the middle of the 1860s on  ``the demand of free mobility  for rigid shapes without change of form in all parts of space''\footnote{``Forderung einer unbedingt freien Beweglichkeit in sich fester Figuren ohne Form\"anderung in allen Theilen des Raumes'' \cite[614]{Helmholtz:1868b}}. 
Twenty years earlier he had reflected on basic features of the concepts of space and time in a manuscript on  ``Allgemeine Naturbegriffe'' (general concepts of nature) which was published posthumously by Leo Koenigsberger \cite{Helmholtz:1845/1903}. In these notes Helmholtz understood {\em space}  as the {\em relationship} for ordering different {\em objects}. In these reflections ``objects'' were understood in a  material sense, and space had to be specified such that ``all kinds of changing  matter distribution'' could still be taken into account.\footnote{``Ferner muss der Begriff des Raumes auch noch so bestimmt werden, dass er alle m\"oglichen Aenderungen der Materie umfassen k\"onne, die hier offenbar nur insoweit in Betracht kommen, als sie Aenderungen der Raumverh\"altnisse, d.i. Bewegungen sind.''
\cite[134]{Helmholtz:1845/1903} \label{Helmholtz 1845}}
 In his later investigations  such changes were restricted to  rigid motions and abstraction was made from the material nature of bodies. Rigidity was now analyzed in terms of geometric (more precisely metric) properties.
 
 He was not the only one to consider rigid  motions as a crucial new feature in the foundations of geometry. Already Riemann had mentioned them among different possibilities for specifying the structure of spaces of constant curvature among all possible (Riemannian) manifolds. Of course he  did not  use the wording of free mobility; he rather  circumscribed the same idea by  postulating  an  ``independence of  spatial figures (K\"orper) of their position'' \cite[283]{Riemann:Hypothesen}. 
 In  France, Ho\"{u}el emphasized the role of rigid motions as an additional foundational feature in \cite{Houel:1867}.\footnote{``Nous demanderons qu'une figure invariable de forme puisse \^{e}tre transport\'ee d'une mani\`ere qulconque dans son  plan ou dans l'espace'', \cite[7]{Houel:1867} here quoted from  \cite[198ff.]{Henke:Diss}.}
That was a year before Helmholtz gave his famous talk at the Heidelberg Medizinische Verein \cite{Helmholtz:1868a}, from which his longer article in the {\em G\"ottinger Nachrichten} resulted \cite{Helmholtz:1868b}.

Peculiar for  Helmholtz's approach was  that the concept of {\em free mobility} stood at the center of his analysis. He established a set of axioms for free  mobility  from which, so he claimed, one could derive the infinitesimal  expression for  distances  by a (non-degenerate positive definite) differential form of second order, similar to what  Riemann had proposed  in the more general approach of his inaugural lecture. In this way he achieved a justification of Riemann's choice for what Weyl would later call the ``Pythagorean nature'' of the metric \cite[\S 4]{Helmholtz:1868b}.
In Helmholtz's analysis the Riemannian metric had to be of constant curvature in order to satisfy the axioms of free mobility. In contrast to Riemann,  Helmholtz  erroneously concluded  from constant curvature and the assumption of infinity that his space had to be Euclidean \cite{Helmholtz:1868a}, \cite[\S 4]{Helmholtz:1868b}. He  was reminded of the possibility of non-Euclidean space structure by  Beltrami in a letter from 24. 04. 1869. The story of this premature conclusion and its correction  has been told at several places;   also the problematics in Helmholtz's statement of his axioms for free mobility.\footnote{\cite{Koenigsberger:Helmholtz,%
Volkert:Helmholtz,Boi_ea:Beltrami,Merker:Espace}.} 

Under the influence of the rising study of transformation groups and their role in geometry by S. Lie and F. Klein, Helmholtz's analysis was received and generalized by mathematicians in different ways. 
One  may even conceive of Klein's idea to characterize a geometry by a manifold and a transformation group ({\em Hauptgruppe}) as  kindred in spirit to Helmholtz's more specific enterprise, even though    Klein did not establish a direct and explicit link to   Helmholtz's analysis of the PoS. He did not even discuss   rigid motions in  his {\em Erlanger Programm}.  We shall see that Klein's view played a crucial role  for E. Cartan's adaptation of the analysis of space to the new situation after the advent of general relativity. 

Sophus Lie, on the other hand, continued  along the lines  opened by Helmholtz. He first sketched  how Helmholtz's  vaguely formulated conditions of free mobility could be restated in terms  of infinitesimal transformation groups \cite{Lie:Helmholtz,Lie/Engel:Raumproblem}; later he investigated it  in more detail  in collaboration with F. Engel   \cite{Lie/Engel:Raumproblem}.\footnote{Cf. \cite{Merker:Espace}.}
 Helmholtz did not  clearly distinguish between  finite and infinitesimal motions. Lie realized that for his goal the interplay of finite and infinitesimal motions played a crucial role. He  
 (and Engel) defined  {\em free mobility}  by a group $G$ operating transitively on a three-dimensional manifold, such that  (1) for any two points $P$ and $P'$ and any two  ``linear elements'' passing through them the group allows to transform these into another, while still ``one continous movement'' is  possible, i.e., one degree of freedom is still open. (2) Only if one fixes two ``surface elements'' passing through the two line elements, the respective transformation is uniquely determined \cite[\S 98]{Lie/Engel:Raumproblem}.\footnote{The modernized description   of a {\em simply transitive} operation on the  {\em flag} manifold is very  close to Lie/Engel's characterization. Note, however, that these authors did { not} demand {\em simple} transitivity, but allowed underdetermination of the operation by the flag up to a discrete subgroup.}
With that more precise definition  and with his general theory of continuous transformation groups as background, Lie (supported by Engel) was able to demonstrate that the space is Riemannian with constant curvature, if free mobility in continuous spaces of dimension $2 \leq n \leq 4$ is assumed.  

Independently Wilhelm Killing started to rethink Helmholtz's analysis under the influence of Weierstrass' approach to mathematics \cite{Killing:1880,Killing:1885}.
 That led him to a mathematical reconstruction of the Riemannian metric from infinitesimal motions  and an  extended study of {\em space forms}, a term  formed by  him, originally with reference to the work of Clifford and   Klein who first took global topological features into accout (which he himself did not).\footnote{For the first aspect see \cite[sect. 4.5]{Hawkins:LieGroups} and \cite{Hawkins:1980};  the second one is discussed by K. Volkert (this volume).} For Killing's peculiar combination of infinitesimal considerations (Hawkins) with those of  the global structure of the geometries  (Volkert) it was a  lucky circumstance that he did   not yet dispose of an explicit concept of transformation groups. He rather characterized the  structure morphisms (to use modernized descriptive terminology) of  space forms  in terms of ``motions of a rigid body''. In hindsight this seems to have came down to considering a groupoid of local transformations containing a neighbourhood of unity of a Lie group and  its infinitesimal transformations 
 (Lie algebra). In the following the historical terminology ``infinitesimal group'' and the present one (``Lie algebra'') will be used 
 indiscrimately. \footnote{The terminology ``Lie algebra'' was introduced by   Weyl in his  lecture course on continuous groups 1933/34;  see \cite[32, fn. 12]{Kosmann:Noether} and, in more detail, the commentary in \cite{Eckes:Commentaire}. \label{Lie algebra}}  
 
 For the next generation, and  already at the life time of Killing, it was Lie's point of view which became dominant. According to Hawkins it was essentially through E. Cartan's work that the part of Killing's contribution relating to infinitesimal transformations entered the mainstream of transformation group studies \cite[290]{Hawkins:1980}.
All in all, the classical space problem had been stated and solved by Riemann, Helmholtz, Lie/Engel, Killing , Clifford, and Klein  between 1851 and 1890. Its outcome showed the following:
\begin{itemize}
	\item  Quite general constraints on the automorphism group of spaces (free mobility in the sense of Helmholtz-Lie/Engel) were able to specify the type of metric (Riemannian) and its structural specification  (constant curvature). The framework and the result of this analysis fitted well to the recent researches on the foundations of geometry.\footnote{For some of the more philosophically inclined participants in the discourse the outcome of the analysis of space seemed to go 
 in hand  with   Leibniz's thoughts on the homogeneity of space. On {\em homogeneity} of space in Leibniz's philosophy see \cite[198ff.]{deRisi:AnalysisSitus}.
} 
	\item More general geometries, respectively space forms, came into sight if one allowed for more general transformation groups as automorphisms  (Klein),  which were no longer bound to physical interpretations,  or for their infinitesimal versions (Killing).
	\item Still other aspects came to the fore, if one took topological characteristics in the large (space forms) into account (Clifford-Klein).
\end{itemize}
  In the decades about the turn of the century  the role of motion groups for the  characterization of  space became widely known in France, most prominently  by talks and popular essays  of  H. Poincar\'e, e.g. \cite{Poincare:1898}.\footnote{Cf. \cite[38--59]{Gray:Poincare}.} His view differed from Riemann's and Helmholtz's by taking sensual and psychological features of the formation of space perception much more serious than the latter, and by his conventionalist  evaluation of the relationship between empirical evidence and conceptual structures.\footnote{\cite{Torretti:Geometry,Walter:2009Poincare}.} 
  
Also  in physics the role of Euclidean motions attracted more and more attention.   In a note written for the third volume (2nd edition) of Appell's {\em Trait\'e de m\'ecanique rationelle}, E. and F. Cosserat demanded to consider action principles  invariant under the group of Euclidean displacements for mechanics, in particular for elasticity theory, and even for (pre-relativistic) physics  in general (``toute la physique th\'eorique'').\footnote{``\ldots on peut observer que l'action, telle que Maupertuis l'a introduite \`a la M\'ecanique, est  {\em  invariante dans le groupe des d\'eplacements euclidiens}. Ce m\^{e}me caract\`ere se retrouve dans la statique des corps d\'eformables, que repose sur la consid\'eration du $ds^2$ de 
l'espace. \ldots Suivant cette id\'ee philosophique, toute la M\'ecanique classique et toute la physique th\'eorique paraissent pouvoir se d\'eduire de la notion unique d'{\em action euclidienne}.'' \cite[557f.]{Cosserat:Note}.  ``Action euclidienne'' was, of course, the  abbreviated expression for an action invariant under the group of Euclidean motion; it stood at the center of the investiations of the Cosserat brothers. \label{fn Appell}}
By this remark they referred back to a proposal made  by Helmholtz. Later it   would come to  play an important role for the theory of generalized deformable  media \cite{Brocato/Chatzis:Cosserat} and for the infinitesimalization of geometric transformation structures in E. Cartan's work  (see section 5).

On the other side of the Rhine, Hilbert  generalized  the classical PoS by passing  over to a topological point of view. In two notes on the foundations of geometry, written in 1902 and 1903, he characterized the (real) plane by  axioms for neighbourhoods, anticipating the later axiomatic characterization of (twodimensional) manifolds, and added two axioms for a purely topological characterization of the group of motions: (I) infinite cardinality of the isotropy group of every point, (II) 3-closedness of the transformation group \cite{Hilbert:GG1902,Hilbert:GG1903,Hilbert:GGAnhangIV}\footnote{See \cite{Strantzalos:Hilbert} (difficult to access). A short description can be found in the historical commentaries in \cite[{\bf 2}, 708ff.]{Hausdorff:Werke}, more in  volume {\bf 6} (in preparation). A transformation group is called {\em 3-closed}, if for any sequence of its operations with images  of a point triple $(A,B,C)$ converging to the point triple $(A',B',C')$, there is  a transformation $g$ of the group, such that $g(A)=A', g(B)=B',g(C)=C'$. }
Hilbert's initiative turned out to be of long-lasting influence on the  topological transformation  group approach to the foundations of geometry  in the course of the 20th century.

The young modernist F. Hausdorff, interested in all aspects of the foundations of geometry, developed his own philosophical outlook on the problem of space trying to keep balance between empirical input in concept formation and mathematical liberty of theoretical generalization (characterized by himself as the view of ``considered empiricism''). In his inaugural lecture  Hausdorff put the discussion of the classical PoS in the wider perspective of general space concepts in the sense of differential geometry and the embryonic theory of point set topology \cite{Hausdorff:Raumproblem}.\footnote{See  \cite[vol. {\bf 6} (to appear)]{Hausdorff:Werke}.}

At that time, the study of electrodynamics and the motion of the recently detected  electrons was already shifting  physics in  a different direction. Hausdorff was well aware of this development. He remarked that ``the hypothesis that the socalled rigid bodies may be subject to certain fine deformations'' had been stated in order to explain the Michelson interference experiments \cite[10]{Hausdorff:Raumproblem}. He even saw the possibility that empirical reasons might speak in favour of  giving up homogeneity and isotropy of space.  
 These were foresightful remarks;  the whole conception of the relationship between mathematical spaces and physical geometry would soon be overturned by the two theories of relativity, special (1905ff.) and general (1915ff.).

\subsubsection*{3. New questions for the space problem raised by relativity}
The rise  of Einstein's theory of special relativity (STR)  led  beyond the classical frame of geometry, although at the outset Einstein did not  intend to do so. In his famous article on {\em electrodynamics of moving bodies} he rather declared coordinate systems to be established ``by means of rigid measuring rods and using the methods of Euclidean geometry'' \cite[892]{Einstein:1905}. 
From Poincar\'e's perspective the situation seemed different. Although he was the first to understand Lorentz transformation in their general mathematical form (in early 1905) his conventionalist philosophy of space led him to believe  that the ``Lorentz-Fitzgerald contraction'' was a real obstacle  for the ``transport of rigid rulers upon which length measurement depends''\cite[204]{Walter:2009Poincare}.\footnote{That this was not the case could be seen most clearly a little later by Minkowki's four-dimensional spacetime (see below).}
Even worse, the discussion of  bodies  in  accelerated motion  showed that   rigidity in a material sense became precarious. The idealization of  putting an extended body into motion instantaneously was, of course, incompatible with the finite speed of light and of causally transmitted  signals more generally; moreover the unequal relative motion of different parts of  bodies, like in the rotating disk, brought  their own difficulties with it. 

 On the other hand, a more abstract mathematical understanding of ``rigid motion'' in the sense of metrical automorphism was only changed, not undermined,  by special relativity, although now the orthogonal group had to be generalized to non-definite signature. That was  most clearly shown  by Minkowski's  conceptual integration of space and time in his four-dimensional  ``world geometry''  \cite{Minkowski:Raumzeit}.\footnote{Cf.  \cite{Walter:2010Minkowski}.}
 The role of  metrical autormorphisms which characterized rigid body motions of the classical,  Euclidean or non-Euclidean, geometries was taken over by  Lorentz transformation (plus translations). Some physicists appreciated this view immediately,  among them  M. Born who considered not only rigid uniform motion but defined special types of rigid accelerated motions by proper time independent ``deformation matrices'' of a material subtrate \cite{Born:1909}.  
 But not all physicists gave due credit to this shift of affairs.  Einstein, e.g., did not realize the conceptual innovative achievement of Minkowski's four-dimensional representation of special relativity before 1912.

 In a talk at London in the summer of the same year,  Poincar\'e, distinguished between the groups of displacement of bodies (in the sense of Helmholtz-Lie) and the group underlying the ``principle of relativity'' for mechanics or electrodynamics. He discusssed the {\em Galilei group}  and the inhomogeneous Lorentz (later {\em Poincar\'e}) { group} as two distinct possibilites. Due to his conventionalist conviction he still argued in favour of the first one for the basic convention regarding spacetime and admitted the second only as a hypothesis of limited range for the study of electrodynamics  and the moving electron \cite{Walter:2009Poincare}, \cite[chap. 6]{Gray:Poincare}. Although from the point of view of the Minkowski/Einsteinians (and in hindsight) this distinction was unjustified, Poincar\'e indicated here in his peculiar way the possibility that  the {\em physical automorphisms} of a dynamical theory  (to take up here a terminology later introduced by Weyl) need not necessarily be identical to  those of spacetime.\footnote{In his later reflections on the problem of space, Weyl posed the question which of the automorphism are of physical import, and which concern the underlying mathematical structure only  \cite[87]{Weyl:PMNEnglish}, a question difficult, perhaps even impossible to decide but  of interest to philosophy of science (see end of section 4). Here, we can consider Poincar\'e's  conventional characterization of the Galilei transformations as geometrical autormorphisms, distinct from the physical transformations of electrodynamics (Lorentz/Poincar\'e group) as a distinction of analogous type. }
 
 Other mathematicians seemed to be more impressed by the basic physical role attributed in special relativity to the Lorentz-Poincar\'e  group. Already the year before Poincar\'e gave his London talk,  G. Herglotz  found that the invariance of the Lagrangian of a relativistically  moving body under the full automorphism group of Minkowski space resulted in dynamical {\em conservation laws}, generalizing those known from classical mechanics:
 \begin{quote} Corresponding to the ten-dimensional group of ``motions'' in $(x,y,z,t)$-space with the metric
 \[ds^2 = dx^2 +dy^2 +dz^2 - dt^2 \]
 10 theorems analogous to the principles of the center of inertia, of the area [conservation of linear momentum, respectively angular momentum, E.S.] and of energy in ordinary mechanics hold for the total body. \cite[511]{Herglotz:1911}$^{(i)}$
 \end{quote} 
 
Although first insights on the relationship between conserved dynamical quantities and invariance under  transformations in Euclidean geometry had been collected during the 19th century (Jacobi, Hamel), it was  the  
 challenge of relativity which led to reconsidering the role of  ``motions'' in classical and in relativistic mechanics.\footnote{In his 1842/43 lectures on mechanics in K\"onigsberg Jacobi derived the principles of ``center of inertia'' and of ``area'' (i.e., conservation of linear and angular momentum) from the observation that the Lagrangian of classical mechanics is invariant under the change of orthogonal {\em coordinate} transformations (shift of the origin and rotation of coordinate axes) \cite[{\bf 2}, 56ff.]{Klein:19Jhdt}. In 1897 I. Sch\"utz derived energy conservation from the dynamical equations of point particle systems. In 1904 G. Hamel  introduced Lie brackets into the study of ``virtual displacements'' in mechanics. Apparently, neither of them recognized energy conservation as a consequence of the time invariance of the Lagrangian \cite[35]{Kosmann:Noether}. } 
   On the advice of  Klein, F.  Engel worked back from Herglotz's relativistic conservation principles to classical ones, by letting the velocity of light $c \rightarrow \infty$. In a second step, documented in a letter written to Klein and published in {G\"ottinger Nachrichten} 1916, he derived energy conservation in Hamiltonian (classical) mechanics from time invariance, avoiding the variational methods of the Lagrangian  approach \cite[36]{Kosmann:Noether}, \cite{Kastrup:1987}.
  
In such a context and in Minkowski spacetime,   Helmholtz-Lie's  flag transitivity, to put it in modernized terminology, lost its interest for researchers as a criterion of free mobility. 
It is unknown to me, whether an attempt to modify  flag transitivity by distinguishing between flag elements with timelike and spacelike markers has been made by any author of the time. Probably nobody did; physicists were interested in dynamical questions, and mathematicians were soon to be  attracted by the problems posed by the generalized theory of relativity (GR). This holds at least for the main authors considered here, H. Weyl and E. Cartan.

General relativity posed, of course, a much greater challenge to the  characterization of the problem of space. Free mobility of finitely extended rigid figures became meaningless in the general case. It could be upheld, if at all, only in the infinitesimal. But even then one could opt for different strategies for posing the question in the new framework. Cartan's main goal was  to integrate Klein's  characterization of spaces (Erlangen program) with Riemann's approach to differential geometry, and to implement it on the level of the infinitesimal structures.
The G\"ottingen mathematicians, including Weyl in his Z\"urich ``outpost'', preferred to stick to the Riemannian approach to differential geometry,   enhanced by Ricci's and Levi-Civita's ``absolute'' differential calculus (i.e, vector and tensor analysis) and, in the case of Weyl, by a generalization to an even more general  ``purely infinitesimal'' founded gauge geometry. 

At G\"ottingen  the question was being discussed,  how the newly acquired insight into the group theoretic derivation of conservation principles (in special relativity and in  classical mechanics) might be adapted to the general relativistic context. Motivated by F. Klein and D. Hilbert, Emmy Noether took up the  question with great success. She found  a general and deep solution to this problem in her now famous publication on {\em Invariant variational problems} \cite{Noether:1918}.\footnote{For an English translation of Noether's 1918 article, a mathematical commentary  and a study of its long delayed reception, see \cite{Kosmann:Noether}.}
Her result consisted essentially of two theorems: the first theorem  generalized the derivation of conservation laws in classical and special relativistic mechanics; it stated the existence of conserved quantities to every generator of a finite dimensional Lie group under which a Lagrangian dynamics remains invariant (and vice versa). Her second theorem dealt with the  infinite dimensional (infinitesimal) symmetries of a Lagrangian. It established differential identies (``Noether identities'') from which, among others, vanishing divergence expressions could be derived. Their interpretation as physically meaningful conserved quantities (currents) was, however,  a much more difficult task than in the  case of special relativity or other fields. It soon led to controversial views. In particular, and most irritatingly for the protagonists, it turned out that energy-momentum conservation could not be generalized to the general relativistic case in  a straightforward 
way\footnote{\cite{Rowe:Noether,Brading:Noether,Brading:2002}.}.

For studying energy-momentum conservation two infinitesimal geometric equivalents to translational invariance have been considered during the last century. In the first one, dominant in the literature on GR until today, infinitesimal diffeomorphisms  were considered.\footnote{Historically, one considered (differentiable) coordinate changes. The ``active'' aspect (diffeomorphism) and the  ``passive'' one (differentiable coordinate change)  can be translated into another and lead to   the same Noether identities.} 
The second  one  developed from Cartan's approach in which point dependent translations are considered  as  part of the structure anyhow. The second view came to maturity only in the second half of the 20th century (see outlook).  The diffeomorphism approach did not lead to an observer independent, covariantly conserved  energy-momentum current; but the corresponding Noether equations (theorem II) are structurally contentful and turned out to be equivalent to the contracted Bianchi identities of the metric.  Hilbert  spoke of ``uneigentliche  Energieerhaltung'' and considered such ``improperness'' as  characteristic of GR \cite{Brading/Ryckman:Hilbert}. 

For fields satisfying the  Einstein equations (in physicists idiom, for fields ``on shell''),  the Noether identities  also imply    vanishing of the  covariant divergence of the energy-momentum tensor, $\nabla_j T^j_k=0$. That was and is a valuable formal property, and could be considered as kind of consolation; but in general  it does not integrate to a conserved energy-like quantity.\footnote{In spite of that, $\nabla_j T^j_k=0$ is often called  ``infinitesimal energy-momentum conservation'', which may appear misleading, if not commented further.}
Proposals for adding additional non-covariant (i.e., non-tensorial) terms $t^i_j$ (energy-momentum {\em pseudotensor}) derived from the gravitational  potentials ($g$ or $\Gamma$) were made by many physicists, among them already quite early Einstein himself \cite{Einstein:Energiesatz}, later among others \cite{Landau/Lifschitz:Fields}. None of these led to a convincing, unanimously accepted, physical interpretation. 

In still another approach   a  so-called {\em canonical} {\em energy momentum} current $\mathfrak t_{i j}$, derived from the  Noether identities, was considered directly. It  is conserved if the fields satisfy the Einstein equations (``on shell''), but it lacks symmetry. Therefore  $\mathfrak t_{ij}$ could  not  be identified, without further ado, with the usual (the  {\em dynamical})  energy-momentum tensor derived from varying  the matter Lagrangian, $T_{ij}=\frac{\delta \mathfrak L_m}{\delta g^{ij}}$. F. Belinfante and L. Rosenfeld  made some progress along these lines.\footnote{
These two authors managed to symmetrize $\mathfrak t_{ij}$  by adding a term $\mathfrak \nabla s_{ij}$ in which   $s_{ij}$  was a linear combination of three  intrinsic (spin) angular momentum terms, $T^{(BR)}_{ij}= \mathfrak t_{ij} + \nabla \mathfrak s_{ij}$. Moreover according to their analysis, the symmetrized tensor  could be identified with the dynamical one,  $T^{(BR)}_{ij}=T_{ij}$  \cite{Rosenfeld:1940}. I thank F. Hehle who did his best to explain these things, well known to any educated relavitst, to me. }
 That was, however, already more than twenty years after the rise of relativity.  So the question whether infinitesimal coordinate changes  (diffeomorphisms)  or differently conceived  infinitesimal transformations ought to be considered as physical morphisms or, rather, as morphisms of the embedding mathematical structure only (to express it again in Weyl's later terminology) was wide open. It  became a point of discussion among philosophically inclined experts in gravity theory.
  \pagebreak

All in all, the rise of relativity, special and general, enriched the analysis of the PoS in many respects. It led to a deeper understanding of the physical symmetries of  classical space (and time) and posed new problems. Some of them remain challenges until today.

\subsubsection*{4. Weyl's analysis of the PoS }
 Weyl turned towards general relativity in 1917,\footnote{\cite{Sigurdsson:Diss,Eckes:2013}.} Cartan a little later.
Weyl chose to take up and to generalize two motifs of the classical PoS: to justifiy the Riemannian nature of the metric (or a generalization of it) in general relativistic spacetimes by group theoretic methods, and to mimic free mobility of rigid matter by allowing what he considered  the ``greatest liberty'' for the distribution of (non-rigid) matter  in space and time. 

In 1918 Weyl developed his {\em purely infinitesimal geometry} which  generalized  Riemannian geometry by introducing a slightly more involved concept of a differential metric, a {\em gauge metric}. Such a gauge metric can be characterized by equivalence classes $[(g,\varphi)]$ of conformally related pseudo-Riemannian metrics $g= (g_{ij})$ and  differential  1-forms $\varphi=(\varphi_i)$, where the equivalence relation was given  by (scaling) {\em gauge transformations} \cite{Weyl:InfGeo}.\footnote{The shortest possible, although only formal, explanation of the equivalence is:  $(\tilde{g},\tilde{\varphi}) \sim (g,\varphi) \longleftrightarrow \tilde{g}= \Omega^2 g \, , \quad \tilde{\varphi}= \varphi - d \log \Omega \, ,$
with $\Omega$ a strictly positive function on the manifold (spacetime). Weyl's  underlying idea was that units of measurement are given by ``local'' (point dependent) specifications, which  allow no immediate comparison. Comparison needs a new structure, the length connection. For more detail see  \cite{Weyl:InfGeo,Weyl:RZM5}, \cite[chap. 5.2]{Adler/Bazin/Schiffer} or   \cite{Scholz:AdP}, embedded in a more historical perspecitve  \cite{Ryckman:Relativity,Scholz:InfGeo}.}
That allowed a direct comparison of lengths of  vectors only if the vectors were attached to infinitesimally close points. For  such a length comparison    between different but infinitesimally close points,  separated by $\delta x$,  Weyl introduced a differential form $\varphi$  which served as an infinitesimal {\em length connection} ($l\mapsto l' = l + \delta l$ with $\delta l \sim - \varphi(\delta x)$ and $\delta l \sim l$). As a crucial ingredient of the geometry, Weyl found that there is {\em exactly one affine connection} $\Gamma$ compatible with the gauge metric. That is,  a vector $v$ of length $l$, infinitesimally parallel translated by $\Gamma$ along the infinitesimal path $u$ (again a vector in the tangent space), changed its length by the amount $\delta l = - \varphi(u) l$.\footnote{$\delta l = |v+ \delta v| - |v|$, where $ \delta v =  - \Gamma(u)v$ or, in index notation, $ \delta v ^i = - \Gamma^i_{jk}u^j v^k$ \cite[12]{Weyl:ARP1923}. With the covariant derivative $\nabla$ with regard to $\Gamma$, compatibility  of $\Gamma$ and the Weylian metric $[(g,\varphi)]$ can be restated by the condition $\nabla g + 2 \varphi \otimes g =0$. For more details see, among others, \cite{Scholz:Connections,Ryckman:Relativity,Eckes:Diss}.}
He called this result the {\em fundamental theorem} of infinitesimal geometry.

That was the basis for Weyl's first attempt of geometrically unifying  gravity and electromagnetism. It  brought him praise and criticism from the side of Einstein and other physicists.\footnote{See \cite{Vizgin:UFT,Goenner:UFT,Goldstein/Ritter:UFT,Scholz:intertwining,Eckes:2013}.}
Einstein praised Weyl's conceptual and mathematical ingenuity, but criticized the latter's attempt to ``improve'' the field theories involved. Among the various points of differences discussed in their correspondence  during the year 1918, Einstein could not see why in Weyl's view  only length standards had to be ``localized'' (present physicists' language), while angle measures were kept immune against the purely infinitesimal relativization. He even ironized that, if one accepted the first, there could just as well come a particularly clever   ``Weyl II'' who might propose localization of angle measurements \cite[551, 31.06.1918]{Einstein:Corr1918}.

 Although  at this time Weyl was probably already convinced  
 that similarity and congruence  at each point of the spacetime manifold had to be specified  by some group theoretic argument, he did not immediately counter this point of Einstein's criticism. But the challenge may  have contributed to his more conceptual analysis of the foundations for  infinitesimal geometry, to which he turned in 1921. At that time  he started to loose his faith in an immediate success for his unified field theory.\footnote{\cite{Scholz:Weyl_matter}}
 
 In the following year, a reedition of Riemann's inuagural lecture gave him the chance to comment on some fundamental aspects of differential geometry  \cite{Riemann/Weyl:Hypothesen}. Riemann had selected the (square root) of a positive definite differential form  $ds^2 = \sum g_{ij}dx_idx_j$ (modernized symbolism) among all homogeneous expressions  $f(\lambda dx_1, \ldots, \lambda  dx_n) = |\lambda |  f(dx_1, \ldots , dx_n)$ ($f$ differentiable, positive), later  called {\em Finsler metrics}.\footnote{\cite{Finsler:Diss}} 
 But he had given  only a quite pragmatic argument for his choice of specialization.

Weyl  conjectured that among the wider class of Finsler metrics the ``Pythagorean'' ones were singled out by the condition that to any such metric  there exists a uniquely determined compatible affine connection, i.e. essentially a parallel displacement of vectors in the manifold.\footnote{Two years earlier Levi-Civita had proved that each Riemannian manifold admits exactly one parallel displacement compatible with the metric, i.e. preserving  lengths of vectors and angles between them.  
Weyl himself had distilled from Levi-Civita's construction  a  symbolic definition of an  abstract ``affine connection'' which  allowed to introduce a  kind of parallel displacement in any differentiable manifold and a kind of differentiation adapted to this geometrical constellation, called a ``covariant derivative'' \cite{Reich:Connection,Bottazzini:Connection,Scholz:Connections}. For a detailed discussion of Weyl's analysis of the space problem in the light of Finsler metrics see \cite{Coleman/Korte:DMV}.} 
It has to be added that  Weyl called a (Finsler type) metric of {\em Pythagorean} nature, if it could be described by a non-degenerate quadratic differential form of any signature (in particular of Lorentz signature). In this way, he transformed the core insight  of what he had identified as the fundamental  of infinitesimal geometry into a criterion for a ``good'' differential geometric metric. He did not pursue this conjecture further in this form,\footnote{The conjecture  (Weyl's ``first space problem'') was proven  40 years later by \cite{Laugwitz:ARP}.  Freudenthal later observed that, in a way, it could be considered solved by the answer to Weyl's (second) PoS \cite{Freudenthal:ARP1960}. In this discussion, he overlooked the specific difference between Weyl's gauge geometry and Riemannian geometry, but the core of his argument is justified:  Weyl's proof also works if one restricts the structure group $G^{\ast}$ of his PoS to the congruence group  $G$, see below.}
but turned towards an even deeper conceptual analysis of metrical structures for infinitesimal geometry. The motif of a {\em uniquely determined affine connection}, however,  stabilized and reappeared also in his more fundamental investigations on the PoS.

In the rewritten fourth edition of {\em Raum, Zeit, Materie} (manuscript finished in November 1920)  Weyl added a new section on the {\em group theoretical outlook on the spatial metric} \cite[\S 18]{Weyl:RZM4}. Here he sketched his new view at the PoS,   introduced  general charachterizations of infinitesimal congruences at any point $p$ and between two infinitesimally close points $p, p'$, and formulated two  postulates P1, P2 (later  postulate of freedom (P1) and of coherence (P2) \cite[lecture 7]{Weyl:ARP1923}). He evaluated the consequences of the postulates for the infinitesimal Lie group  (Lie algebra) $\mathfrak g$ and conjectured that the only groups satisfying the postulates might be  the generalized special orthogonal groups $\mathfrak g = so (p,q)$  ($p+q=n$ dimension of the manifold) \cite[132]{Weyl:RZM4}.\footnote{Weyl stated that he had been able to prove the conjecture for the cases $n=2,3$, but not yet in general \cite[133]{Weyl:RZM4}}
Optimistically he declared that, if the general proof could be given, ``the essence of space (Wesen des Raumes) would have been made comprehensible''  by deep mathematical considerations. After an eloge on the long historical course of events leading to this result (Euclid -- Newton -- Gauss -- Riemann -- Einstein) he declared in an often quoted remark
\begin{quote} The example of space is most instructive for that question of phenomenology that seems to me particularly decisive: to what extent the delimitation of the essentialities (Wesenheiten) rising up to consciousness express a characteristic structure of the domain of the given itself, and to what extent mere conventions participate in it. \cite[133f.]{Weyl:RZM4} 
$^{(ii)}$
\end{quote}
 
This remark is often taken as proof for Weyl's close affiliation to Husserl's phenomenology.\footnote{For example in \cite[157f.]{Ryckman:Relativity}, from which I have taken the translation of the Weyl quote.}
 But one may doubt that this was the driving force of his thought. Declaring to understand the specific relationship between the {\em structure of the given itself} and the {\em conventional} as ``particularly decisive''   sounded rather like a critical commentary to Poincar\'e's latest discussion of this question in \cite{Poincare:1912}, even though it was stated in a language which clearly did not hide Weyl's involvement in the discourse of German idealistic philosophy. His  epistemic orientation followed a  
 road passing  equidistantly  between  an  empirically oriented positivism and idealism.\footnote{The idealistic aspect of Weyl's view has been emohasized by \cite{Ryckman:Relativity,Bernard:Idealism}.}
  It was  governed by a  hope  expressed by Riemann half a century earlier that    scientific knowledge is able, step by step, to go  ``behind the surface of appearances''. 

In the next two years Weyl elaborated the consequence of his analysis of the PoS \cite{Weyl:ARP1921,Weyl:ARP1922a,Weyl:ARP1922b}. In April 1922, in lectures at Barcelona (in French) and Madrid (partially in Spanish), he had the chance to present his analysis  to a Spanish audience. He edited his lecture notes  in German  as a booklet,  complemented by 12 mathematical appendices  \cite{Weyl:ARP1923}.\footnote{For a detailed discussion of the Barcelona lectures see \cite{Bernard:Barcelona}.}
 This publication  contains the most extended presentation of his thoughts on the PoS and its embedding in the wider mathematical context.

Already in his first presentation of the new PoS  Weyl declared that  one could perceive as many different metrical determinations as there are ``essentially different'' (i.e., non conjugate) subgroups of the general linear group $G \subset GL(n,\R)$. He  posed the question which of these could justifiably be considered as {\em rotations} in a generalized sense \cite[125]{Weyl:RZM4}. Such ``rotations'' had to be considered in the infinitesimal neighbourhoods of points; they operated on the ``vector body (Vektork\"orper)'' attached to each point $p$ separately (today, on the tangent spaces $T_pM$). 
As first constraint he introduced volume invariance, i.e., $G \subset SL(n,\R)$.

If, with regard to a (generalized) rotation group $G$, conjugation by a linear transformation $u$  does not modify the group,  $u^{-1}\cdot G\cdot u = G$, the transformation leaves the ``metrical determination'' correlated with $G$ unchanged, i.e., it could serve as a (generalized)  { similarity}. Therefore Weyl considered the normalizer $G^{\ast}$ of $G$   as the {\em similarity group} of the rotations $G$. A {\em metric connection} between any two infinitesimally close points would be given, so Weyl, by specifying a linear connection $\Lambda$ expressed, after choice of coordinates and of a basis of the  ``vector body'' (tangent space in our terminology) by a  system of real numbers $(\Lambda^i_{jk})$ ($1 \leq i,j\leq n$), not necessarily symmetric in the lower indices (in Cartan's terminology possibly with torsion), and transforming correctly. 

This ``metric connection'' should not be confounded with a parallel displacement  (affine connection) which only stood in indirect relation to it. The attribute ``metric''  only  indicated that a transfer of vectors or figures defined by vectors led to vectors or figures which ought to be considered ``congruent'' in the geometry under consideration. Therefore any composition of a {\em congruence displacement} with a {\em rotation} in the target space  had to be considered just as well as a metric or congruent discplacement. Thus Weyl considered metric displacements to be defined by a linear connection $\Lambda$ up to additive modifications by a point dependent infinitesimal rotation, i.e., up to addition of differential one form $A =(A^i_{ik})$ with values in $\mathfrak g$.\footnote{That means $A = A^i_{jk} dx^k$, $(A^i_{j}) \in \mathfrak g$, using Einstein sum convention. }
In consequence the metric connection was given only up to equivalence
\[ \Lambda^i_{jk} \sim \tilde{\Lambda}^i_{jk} \longleftrightarrow 
\tilde{\Lambda}^i_{jk} = \Lambda^i_{jk} + A^i_{jk} \quad 
\] 
for some differential form $ A^i_{jk} dx^k $ with values in $\mathfrak g$.

All this, including the explanation of the concept of {\em parallel displacement} by an affine connection $\Gamma = (\Gamma^i_{jk})$ symmetric in the lower indices $j,k$, was considered by Weyl as  ``mere conceptual analysis, an explication of what is contained in the concepts metric, metric connection, and parallel displacement as such''. The specific postulates which ought to be satisfied added something more to the analysis, it belonged to ``synthetic part [of the analysis of PoS] in the Kantian sense'' \cite[49]{Weyl:ARP1923}. Roughly formulated the postulates were:\footnote{For more details and Weyl's motivational arguments see, e.g. \cite{Scholz:Weyl_PoS}.}

\begin{itemize}
	\item [P1]  (``Principle of freedom''): 
$G$ allows the ``widest concievable range of possible congruence transfers'' in one point. \ldots
\item[P2] (``Principle of coherence''): 
To each congruent transfer $ \Lambda  ^i_{jk}$ exists {\em  exactly one equivalent  affine} connection. 
\end{itemize}
Weyl explained the meaning of ``widest conceivable range'' 
 by the possibility that the group $G$ had to be large enough that at a given point $p$ the coefficients of the metric connection $(\Lambda^i_{jk})$ could take any  values in the reals, i.e., the connection had degree of freedom $n^3$. 
 
 In the  Barcelona lectures he added that this postulate assured greatest possible liberty for matter to shape the metric and ought therefore  be  considered as the relativistic equivalent to free mobility in the classical PoS. This argument was surprisingly close to Helmholtz's {\em early} (1845) characterization of the concept of space. There Helmholtz had stipulated that the concept of space ought to allow ``all possible changes of matter''. Only twenty later he boilt down this rather general description to free mobility.\footnote{See above fn. (\ref{Helmholtz 1845}).}
 
Because homogeneity of space was the result of free mobility, and free mobility  the expression of  allowing ``all possible changes of matter'' in the classical PoS, Weyl had good reasons to argue:  
 \begin{quote} The possibility {\em to subject the metric field to arbitrary virtual changes, keeping the nature of the metric fixed,} has taken the place of homogeneity of the metric field postulated by Helmholtz. \cite[46]{Weyl:ARP1923}$^{(iii)}$
 \end{quote}
The relation of Weyl's postulate P1 to Helmholtz's free mobility becomes much clearer in the light of the {\em Nachlass} notes from 1845; but it has to be added that we do not have any  information whether Weyl had read Koenigsberger's edition of the  early Helmholtz text. That may have been  the case, but he could just as well have arrived at his generalization by his own conceptual analysis.
 
The two postulates led to constraints for the infinitesimal Lie group $\mathfrak g$, expressed by

\begin{lemma} If a group $G \subset SL(n,\R)$ with infinitesimal rotations (Lie algebra) $\mathfrak g$ characterizes the ``rotations'' in an $n$-dimensional  generalized metric structure, in the sense of the postulates P1, P2,    the following holds:
\begin{enumerate}
	\item $N:= dim\, \mathfrak g = \frac{n}{2}(n-1)$,
	\item $ tr\,  A = 0$ for any $A \in \mathfrak g$,
	\item the only system of $N$ matrices $A_{k}= (A^i_{j k})$ in $\mathfrak g$ ($1 \leq k \leq N$) which is symmetric in the lower two indices, $A^i_{j k}=A^i_{ k j}$, is  identical $0$. \cite[131f.]{Weyl:RZM4}, \cite[51]{Weyl:ARP1923}.
\end{enumerate}
\end{lemma}

In 1920 he conjectured that the only groups satisfying the conditions of Lemma 1 are the  special orthogonal groups of arbitrary signature, $G=SO(p,q, \R)$. 
He was able to give a proof at first in the dimensions $n=2,3$,  then, in a case by case study, for general $n$  \cite{Weyl:ARP1922a},  \cite[app. 12]{Weyl:ARP1923}. 
Already his first case by case study brought him closer to the theory of infinitesimal Lie groups, but according to Hawkins he did not yet delve deeply into it \cite{Hawkins:LieGroups}. 

 Cartan, on the other hand, read  Weyl's first presentation of  the PoS in 1921 or early 1922. He tried to make sense of it in terms of his newly developed ``non-holonomous'' spaces and 
reformulated Weyl's condition in such a way that his classification of infinitesimal Lie groups became applicable. By  this move he could solve a  (modified)  main theorem of the space problem without much additional technical work ({\em Cartan's space problem}, see below). In October 1922 he wrote Weyl about this result. Weyl was impressed and was attracted to Cartan's theory, but did not yet try to assimilate it \cite[438ff]{Hawkins:LieGroups}. He continued to rework  the proof of the extended Barcelona lectures in his own approach. 

In a systematic study of the (complex) rotation group Weyl  brought it into a form which he found finally satisfying \cite{Weyl:1923Drehungsgruppe}.  The result can be stated as
\begin{theorem}
The only complex infinitesimal Lie groups $\mathfrak g$ satisfying conditions (1) to (3) of Lemma 1 are the Lie algebras of special orthogonal groups, $\mathfrak g = so(n,\C)$ ($n\geq 2$). 
\end{theorem}

Over the complex numbers any nondegenerate quadratic form can be brought into normal form with positive coefficients, even normalized to 1, while for real subspaces positive definiteness may be lost. The result was exactly what Weyl had conjectured from the outset. He drew the resum\'ee:
\begin{quote}
Just as the old conception of a metrical structure inbuilt into space by a priori reasons independently of matter leads to Helmholtz-Lie's characterization of the rotation group, the modern view due to Einstein's relativity theory, according to which the metric structure is variable and is causally dependent on matter,  leads to acknowledging a different property of the rotation group as decisive (\ldots). \cite[348]{Weyl:1923Drehungsgruppe}$^{(iv)}$
\end{quote}

For Weyl, the conditions (1) to (3) of the lemma took over a role  exactly analogous  to flag transitivity in the old one.\footnote{Already early in his proof attempts Weyl realized that the trace condition (2) could be omitted for $n>2$ \cite[112]{Weyl:ARP1923}, \cite[354]{Weyl:1923Drehungsgruppe}}
 By complexifiation they   led  to the inclusion of non-definite rotation groups in the real case. 

During the long path toward a satisfying proof of his main theorem Weyl found sufficient occasion to rethink his position with regard to the foundations of mathematics. In 1921 he was an ardent defender of what he considered Brouwerian intuitionism, among others he  considered the principle of excluded middle as illegitimate for infinite collections.\footnote{In fact, Weyl's perception of intuitionism differed from Brouwer's in certain respects \cite{Hesseling:Gnomes,vanDalen:Hendricks}.}
Already in his first systematic case by case proof Weyl realized how difficult proofs could become if one was to avoid the true/false dichotomy of classical logic. As a result he chose to weaken the intuitionistic proof criteria for the continuum in order to arrive at clean case distinctions.\footnote{He declared his own proof as ``\ldots nur zwingend, wenn wir in einem Zahlbereich bewegen, wie z.B. den rationalen Zahlen, in welchem eine Zahl entweder $=0$ oder $\neq0$ ist. Die Fallunterscheidungen der Elementarteilertheorie sind von diesem Standpunkt aus ein besonders bedenklicher Ausgangspunkt.'' \cite[295]{Weyl:ARP1922a} } 
But he still promised to ``come back to this point, once a new analysis in more definite form has been elaborated than is the case at the moment'' \cite[295]{Weyl:ARP1922a}. 

During the work process in the following year he seems to have changed his mind. He was obviously satisfied with his final proof in \cite{Weyl:1923Drehungsgruppe} which still needed dichotomic case distinctions. The only remnant of his formerly declared goal of intuitistic proofs was hidden under a weakly self-ironic remark at the end of the article:
\begin{quote} Our game on the chess board of matrix schemata has been played to the end \ldots (ibid., 372) 
\end{quote}
Readers who knew his commentaries on the foundations of mathematics could easily recognize the parallel to earlier critical remarks of Weyl on Hilbert's formalism as an epistemically empty game (chess). But in 1922/23 Weyl was sure that the game had an epistemical surplus value provided by the context of the problem of space. 

For him, the context went far beyond the algebraic characterization of the rotation groups. With the geometrical characterization of infinitesimal congruence structures by the postulates P1 and P2  it seemed  clear that over the reals the congruences in the infinitesimal were characterized  by a rotation group $G=SO(p,q,\R)$  of fixed signature (by reasons of continuity). After the choice of a basis in the ``vector bodies'' (tangent spaces),  they  could be specified  as those of  a pseudo-Riemannian metric $g= (g_{ij})$. Pointwise change of the congruence groups was mediated by conjugation in the similarity group $G^{\ast} = \R^+ \times SO(p,q,\R) $ (in modernized terminology the structure group). Thus displacement of vectors or a vecor constellation (a figure) by the ``metrical connection'' from a point $p$ to an infinitesimally close one $p'$ would lead to a {\em similar} constellation (figure). The deviation of the parallel discplacement (unique among the all equivalent metrical connections) from conserving metrical values relative to the chosen $g=(g_{ij})$  was then expressible by a Weylian length connection $\varphi = (\varphi_i)$ depending on $g$. 
 Thus Weyl  arrived at
 \begin{result}
The structure of a Weylian metric $[(g,\varphi)]$ can be  reconstructed from  infinitesimal congruences in the sense of the new PoS (postulates P1 and P2), just as the Riemannian metric had been reconstructed by Helmholtz, Lie, and Killing in the classical PoS \cite[51f.]{Weyl:ARP1923}.
\end{result}    
 
 In this sense, Einstein's ``Weyl II'' objection was shown to be immaterial, if one accepted Weyl conceptual analysis and his ``synthetic'' postulates (at least P2).\footnote{In the later language of fibre bundles, Einstein's objection can be stated as the question, why one would not expect a  general affine connection, not reducible to an orthogonal group, if one was prepared to give up Riemannian geometry anyhow. (Physicists speak of a {\em metric affine} theory of gravity.) Eddington's affine field theory went in that direction, even in a more radical form; no wonder that it played a crucial role in Einstein's early turn toward unified field theories.}
 The structure of  Weylian metric had been given a conceptually outstanding status,
  although  Einstein probably would have been  unimpressed by such a  combination of philosophical and mathematical arguments in Weyl's analysis, in case he came to an understanding of the complicated argument.\footnote{A side remark in Einstein's letter from 06. 06. 1922 to Weyl shows that he was aware of Weyl's analysis of the space problem: ``Ich bem\"uhe mich gerade um das Verst\"andnis Ihrer Arbeit \"uber die mathematische Vorzugsstellung der quadratischen Form. Physikalisch komme ich nicht weiter. Ich glaube nicht an den Zusammenhang zwichen elektrischem Feld und Streckenkr\"ummung. \ldots'' \cite[{\bf 13}, doc 219]{Einstein:Papers}. The second phrase shows that Einstein still struggled with Weyl's unified field theory. --- Thanks to  D. Lehmkuhl whom I owe the information about  Einstein's (at least marginal) interest in Weyl's PoS.} 
  At the time of Weyl's work on the PoS Einstein was already on his path toward metric-affine theories \cite{Sauer:EinsteinUFT}, although in some  moments he still tried to make sense of Weyl's field unification.\footnote{In his travel diary Oct. 1922 -- March 1923 (Japan, Palestine, Spain) he noted  (09. 10. 1922): ``Gestern habe ich die elektromagnetischen Vakuum-Gleichungen (\ldots) im Sinne der Weylschen Geometrie umgerechnet, in der Hoffnung, einen Ausdruck f\"ur die Stromdichte zu finden. Es kommt aber unbrauchbares Resultat $\varphi^{\mu \alpha}\varphi_{}\alpha$ heraus.'' \cite[{\bf 13}, doc 379]{Einstein:Papers}. }

  For Weyl, on the other hand, his way of treating the PoS may have contributed to pave the way for his later work in generalizing gauge structures from the similarities to phase factors in Dirac theory (group extension by $U(1)$ rather than by $\R^+$ \cite{Weyl:1929Dirac}.\footnote{\cite{Vizgin:UFT,ORaif/Straumann,Scholz:Weyl_PoS}} 
  Of course this shift from  similarities  to  phase was part of a  much wider turn from classical unified field theories to the question of how to make the views of general relativity coherent with the new wave functions  (Schr\"odinger and/or Dirac) in quantum theory. For Weyl, the turn went in hand with  a rising awareness  that a priori type arguments were only of limited help for physics, while the empirical input had to be considered as the finally decisive factor.\footnote{\cite{Eckes:2013,Scholz:Fock/Weyl,Scholz:Weyl_realism}.} 
It may have been this shift which induced  Weyl to ponder about  the difference between purely {\em geometrical} (in the sense of mathematical geometry) and {\em physical automorphisms} in physical theories.

  In the English translation of his {\em Philosophy of Mathematics and Natural Sciences} and in the unpublished notes of a talk, given apparently about the time of preparing \cite{Weyl:PMNEnglish},\footnote{H. Weyl ``Similarity and congruenc: a chapter in the epistemology of science'', Bibliothek ETH Z\"urich, Hs 9a:31 (23 pp.), formally undated typoscript,  written ``about 30 years'' after ``an attempt was made by the speaker to build up a unified field theory of gravitation and electromagnetism \ldots'' (ibid, p. 20).  }
   Weyl  clearly distinguished between different types of symmetries in physical theories:  between those that are of mathematical, in particular geometrical, importance ({\em mathematical} or {\em geometrical autormorphisms}) and those that relate to the ``physical world'' ({\em physical automorphisms}). 
  
  Discussing  classical geometry and classical physics he explained the group  of Euclidean isometries, introduced as proper or improper ``motions'', and  denoted it by $\Delta$. He then went on:
  \begin{quote} A far deeper aspect of the group $\Delta$ than that of describing the mobility of rigid bodies is revealed by its role as the group of {\em automorphisms of the physical world}. In physics we have to consider not only points but also various other types of physical quantities, velocity, force, electromagnetic field strenght, etc. \ldots All the laws of nature are invariant under the transformation thus induced by the group $\Delta$. \cite[82f., emphasis in original]{Weyl:PMNEnglish}
  \end{quote}
On the other hand, he  introduced the   (classical) similarity group $\Gamma$ as the  normalizer of $\Delta$ (in the diffeomorphism group of space) and explained  its importance for characterizing  classical metrical (Euclidean) geometry. Different to what he had thought at the turn to the 1920s, he now drew upon the  empirical insights  from atomic physics:
  \begin{quote}
  The atomic constants of charge and mass of the electron and Planck's quantum of action $h$ fix an absolute standard of length, that through the wave lengths of spectral lines is also made available for practical measurements. \ldots  The orthogonal transformations of signature $-$  must be included in $\Delta$. For there is no indication in the laws of nature of an intrinsic difference between left and right. (ibid.)
  \end{quote} 
 In this sense the group $\Delta$ characterized ``physical automorphisms''. For the geometrical (mathematical) automorphisms Weyl concluded now:
  \begin{quote} The group of geometrical automorphisms, by virtue of the very meaning of this term, is the normalizer $\Gamma$ of $\Delta$. It is larger than $\Delta$ as it includes the dilatations \ldots. This divergence between $\Delta$ and $\Gamma$ proves conclusively that {\em physics can never be reduced to geometry} as Descartes had hoped [and Weyl himself had believed  between 1918 and ca. 1923, E.S.]. \cite[83]{Weyl:PMNEnglish}
  \end{quote}
  
  In other words, it was the {\em distinction between the automorphism groups} of mathematics/geometry and physics which led Weyl to the  conclusion that physics could not be ``reduced'' to mathematics.\footnote{Although this insight can be defended, and probably would better be  defended, by  different arguments and even may  seem obvious to unbiased reflection, this conclusion was important for Weyl, because in his radical years 1918--1920 he had tended to believe in such reducibility. Weyl's clearcut distinction between physical and mathematical automorphisms became blurred again after the rise of modern gauge theories in late 20th century physics (see below). }
 
 A little later in the chapter Weyl discussed the new situation after the advent of general relativity. He emphasized that the metric is here no longer part of the spatial structure as such (``the extensive medium of the world''), but rather part of its  changing material content. By the argument that  ``the free mobility of bodies without changes in measure is regained, since a body in motion will `take along' the metric field that is generated or deformed by it''  (ibid, 87), Weyl argued for accepting diffeomorphism with dragging of physical fields as part of the physical automorphisms of GR. Moreover, he presentated  the  method of moving (orthogonal) frames,  denoted by $e_{i \beta}, (1 \leq i, \beta\leq n )$, used by Cartan and by himself in his second (phase) gauge theory of the general relativistic Dirac field  \cite{Weyl:1929Dirac}.  Then he went on:
 \begin{quote}
 What has happened in the transition from special to general theory is obviously this. The physical automorphisms forming the group $\Delta$ as described in the previous section\footnote{Weyl skipped the discussion of special relativity as an  intermediate step  and referred the reader to another chapter of his book, if he or she wanted to see `` how time is included as a fourth coordinate in the above scheme'' \cite[83]{Weyl:PMNEnglish}. }
 have been split into their translatory and rotatory parts. The group of translations has been replaced by that of all possible transformations of the coordinates, whereas the rotations have remained Euclidean rotations but are now tied to a center $P$ and must be allowed to vary freely while the center $P$ moves over the manifold. Space, the extensive medium of the material world, is clearly the seat of the group $\Omega$ of coordinate transformations; but the group $\Delta_0$ [the point dependent representation of the rotation group, E.S.] seems to have its origin in the ultimate elementary particles of matter. The quantities  $e_{i \beta}$ thus mediate between matter and space. \cite[89]{Weyl:PMNEnglish}
  \end{quote} 
  
 About 1928/1929,  Weyl's caracterization of the problem of space had changed  from his earlier strongly aprioristic philosophical view  to an   empirically imbued one.\footnote{In the decade 1920-1929 Weyl thus moved from a heavy leaning toward the a priori side of a  ``relativized a priori'' in the sense of \cite{Friedman:Dynamics} to the  historically relativized one, taking account of new empirical constraints (atomic and quantum physics) and taking the empirical input of knowledge more serious than before. }  
  The old gauge idea (scale geometry) now appeared to him  as part of the ``geometrical automorphisms'', not of the physical ones.\footnote{The reappearance of scale in- and covariance in field physics in the second half of the 20th century brought new aspects into this question. These cannot be discussed here.} 
  The physical automorphisms he had now in sight included, slightly rephrased, the diffeomorphism auf the spacetime manifold (``all possible coordinate transformation''),  point dependent Lorentz rotations, extended by (also point dependent) phase transformations, not unlike the similarity extension of 1918 and the early 1920s. Clearly,  the extension group, $U(1)$, was no longer part of the geometry (it did not operate on the ``external medium of the world''),  but operated on the dynamical space of Dirac spinors or the wave fields. This did not hinder Weyl to consider  physical automorphisms with point dependent group operations (i.e., gauge groups of physical theories) as carrying the ``character of general relativity'', \cite[246]{Weyl:1929Dirac}.\footnote{About a decade later, Weyl reconsidered the Lagrangian of the general relativistic Dirac field and found a ``slight discordance between affine connection and metic'', i.e., a (con)torsion effect resulting from the spin current of the Dirac field \cite{Weyl:1950}. He thus realized that there may be physical reasons to relax the condition of {\em affine} connection in general relativity, which in  1923 he had considered a ``synthetical'' postulate ``a priori'' (see section 4). But there is no indication that he started to rethink the question of ``physical automorphisms'' in full-fledged Cartan geometric terms.} 
  
  The tendency of intertwining geometry with matter in Einstein's theory became strengthened with the new, although still embryonic, insights in the physics of elementary particles at the end of the 1940s.   At the time of rewriting his {\em Philosophie der Mathematik und Naturwissenschaften}  Weyl even dared to speculate that the  (Lorentz) rotations and their extension  ``have their origin in the ultimate elementary particles of matter''. Mathematically, the localized group ``could be considered as an abstract group of which various representations by linear transformations are characteristic for various physical quantities'' \cite[89]{Weyl:PMNEnglish}. Groups and their representations, gauge structures, geometry and quantum physics intermingled in Weyl's view  at the turn to the 1950s. The problem of space could no longer be separated from the problem of matter.

 \subsubsection*{5. Cartan's  view of the space problem}
When Cartan got to know  Weyl's  PoS (in 1921 or in early 1922)
he  could already build upon a huge expertise in  the theory of { infinitesimal Lie groups} 
 which he had collected over a period of roughly thirty years.\footnote{In his first note on Weyl's PoS \citeasnoun{Cartan:1922PoS}  referred to the French translation of Weyl's {\em Raum - Zeit - Materie} \cite{Weyl:RZMfranz}. It may very well be that he knew from Weyl's PoS already from the German fourth edition, from which the translation was made; at least he cited the {\em German} title of the book.}
  Among others  he   had classified the simple complex Lie groups in his doctoral thesis 1894,  twenty years later (1913/14) the real ones \cite[part III]{Hawkins:LieGroups}. Moreover he had  brought to perfection the usage of    {\em differential forms} (``Pfaffian forms'')   in  differential geometry  \cite{Katz:DiffForms}. In 1910 he had started to describe the differential geometry of  classical motions by generalizing Darboux' method of
  ``tri\`edres mobiles'' (moving frames)  \cite{Cartan:1910}.\
  During the year 1921 he started to develop, and in 1922 he published, a series of articles in which he presented  his new theory of generalized {\em non-holonomous} spaces (later {\em Cartan spaces}), a kind of ``infinitesimalized'' Kleinian geometries which allowed a supplementary  translational curvature called {\em torsion}.\footnote{See \cite{Nabonnand:Cartan_2009} and P. Nabonnand's contribution on Cartan in  this volume. For the mathematical background see, e.g., \cite{Sharpe:DiffGeo}. }
  
  The idea for such an additional aspect of curvature seems to have been a result of Cartan's attempts to geometrize certain aspects  of E. and F. Cosserat's theory of  generalized   elastic media. The brothers had  investigated a  variational principle for elastic media with an action density dependent on an infinitesimal ``tri\`edre trirectangulaire'' (orthonormal frame) and  invariant under infinitesimal Euclidean motions.\footnote{Compare  fn. (\ref{fn Appell}) and, for more details, \cite{Brocato/Chatzis:Cosserat}. For a more systematic physical discussion see, among  others, \cite{Hehl:2007}.}  
  Under such assumptions the ensuing stress tensor was no longer necessarily symmetric. Moreover,  an ``infinitesimal'' surface element on the boundary of an elastic body under external influences could be the carrier of an elastic torque in addition to an elastic force.\footnote{``\ldots nous appellerons {\em effort ext\'erieur et moment de d\'eformation ext\'erieur au point} $M$ {\em de la surface} $S$, qui limite le milieu d\'eform\'e, {\em rapport\'e \`a l'unit\'e d'aire de surface} $S_o$ les segments issus de ce point $M$ et dont les projections sur les axes $Mx',My',Mz'$ sont respectivement $-F_o',-G_o',-H_o'$ et $-I_o',-J_o',-K_o'$ '' \cite[598]{Cosserat}. The quantities $F_o' \ldots $ and $I_o'$ arose as  coefficients of infinitesimal deformation vector, respectively infinitesimal rotation vector, in the surface action density (after application of Stokes theorem) \cite[597]{Cosserat}.} 
  
  At the end of his {\em Comptes Rendus}  note presented on February 27,  \citeasnoun{Cartan:1922[58]} remarked that this generalization was close to the Cosserat's investigation in \cite{Cosserat:Note} and to Weyl's  generalized concept of space in the PoS.\footnote{``J'ajouterai que les consid\`erations pr\'ec\'edentes qui, du point de vue m\'ecanique,
s'apparentent aux beaux travaux de MM. E. et F. Cosserat sur
l'action euclidienne, s'apparentent \'egalement à la th\'eorie des espaces g\'en\'eralis\'es
de H. Weyl et peuvent elles-m\^emes se g\'en\'eraliser.'' \cite[595]{Cartan:1922[58]}. }
In the light of Weyl's  recent studies of spaces with  infinitesimal congruence structures the conceptual problems of    
general relativity presented a challenging context for developping Cartan's ideas on connections with a translational component. Why should not, perhaps, the relativistic ``ether'' (to use Einstein's and Weyl's expression for the structure field of general relativity) behave   as sophisticated as the Cosserats' solid state media could (at least hypothetically)?
  
As a consequence, Cartan's generalization went far beyond the one envisaged  by Weyl. He developped it in a multitude of publications extending over years, and even decades, to come. As  a side activity to his main work he considered it worthwhile to give his own answer to the question  posed by Weyl the new PoS \cite{Cartan:1922PoS,Cartan:1923PoS}. But his approach differed in three respects from Weyl's:
\begin{itemize}
	\item He reformulated Weyl's characterization of the  PoS in terms of his  representation of connections by differential forms.
	\item As a result of his reformulation he unknowingly suppressed  the specific Weylian aspect of the metric and reduced the geometrical analysis of the PoS to the (pseudo-) Riemannian special case.
		\item The conditions arising from the PoS for the infinitesimal groups then acquired such a form that Cartan's general structure theory of simple and semi-simple Lie groups could be applied, leading to a much shorter proof of the main theorem.
		
\end{itemize}

In a letter to Weyl, written in October 1922, Cartan politely announced his second publication on the PoS and sent him a preprint (``un exemplaire d'un m\'emoire que je viens de faire para\^itre dans le {\em Journal de Math\'ematiques}'') for information.   
  Cartan modestly played down the achievements of his general structure theory of infinitesimal Lie groups,
  but tried to interest Weyl for his methods to describe connections by systems of differential forms which we today prefer to read as differential forms with  with values in a Lie algebra.\footnote{``\ldots une foi le th\'eor\`eme reconnu vrai, le plus ou moin de simplicit\'e de sa d\'emonstration n'est rien aupr\`es de sa port\'ee philosophique \dots    Les proc\'ed\'es de calcul qu j'y emploie ne sont pas ceux du calcul differentiel absolu; \ldots  ils sont tout aussi bien adapt\'es je crois \`a l'\'etude des multiplicit\'es, non seulement \`a connexion affine (pour employer votre terminologie), mais encore \`a connexion projective, ou conforme, etc. \ldots'' \cite[Cartan to Weyl, 09. 10. 1922]{Weyl/Cartan:Korr1920s}. I thank P. Nabonnand for having made transcripts of this letter available to me.}
  
 Weyl's characterization of the PoS, with its  basic ingredient of a  generalized ``group of rotations'' $G$ operating on the vectors attached to a point  (``vecteur issu de $P$'' in Cartan's language) was easy to accept for Cartan. But Weyl's way of dealing with similarity operations   $G^{\ast}$ for the  displacements of vectors  between infinitesimally close points $p$ and $p'$ remained incomprehensible to  him, because Weyl used a  form of group extension ($G \subset G^{\ast}$) specific for his conception of gauge geometries.\footnote{Cartan's  own approach to non-holonomous spaces also used  a group and a subgroup ($L \supset G$), but for the completely different purpose of constructing homogeneous spaces $L/G$ which were to describe infinitesimal tangent structures. See below or, in a systematic perspective, \cite{Sharpe:DiffGeo}.}
 Moreover, the characterization of the postulate of freedom (P1) with the specification that the ``metric connection'' $\Lambda^i_{jk}$ should have ``greatest liberty''  for matter influence and  could acquire any of $n^3$ values (at a point), could not easily be rephrased (if at all) in Cartan's terminology of moving frames and differential forms.\footnote{If one wanted to express Weyl's P1 in modernized terminology, one had to exploit the interplay of a connection with values in the general linear group, which can acquire any value at a point, although it can be reduced to $G^{\ast}$, respectively its Lie algebra $\mathfrak g^{\ast}$. Moreover one has to observe that the   equivalence classes of postulate P2 are constructed by differential forms with values in $\mathfrak g = Lie\, G$. }
 So it is not surprising that after rephrasing Weyl's first postulate,
 Cartan continued with an interpretation of his own.
   
Cartan used  a base choice of the vector spaces ($T_q M$, $q\in U$) in a (finite) neighbourhood $U$ of $p$, such that the operation of the ``rotation group'', or its infinitesimal version, was given by constant coefficients. In other words,  his way of making precise Weyl's speech of ``different orientations'' of the group $G$ at different points consisted in assuming a choice of  frames (vector bases) in each $T_qM$, chosen  such that the operation of $G$ has the same matrix reprentation for all $q \in U$.\footnote{In modernized terminology one could see in Cartan's description a well chosen trivialization of the tangent bundle,  adapted to the operation of $G$ on the tangent spaces.}
  The next step contained Cartan's simplification of Weyl's approach:
 \begin{quote} This given, the meaning of Mr. H. Weyl's first axiom is the following: We may choose at each point $P$ of the space the {\em orientation} of the group arbitrarily, i.e., the coefficients $a_{ij}$  \ldots [which describe the operation of the group, ES]. After this choice, the metrical connection of the space is determined: one of the congruence transfers \ldots is the one in which two corresponding vectors have {\em the same} components \ldots
  \cite[172]{Cartan:1923PoS}$^{(vi)}$
 \end{quote}
Thus Cartan assumed that the frames chosen for the representation of the $G$-operation could be considered as {\em congruent} without further specification. Put in the language applied by him and Einstein in the late 1920s,  he considered a (Cartan-) connection corresponding to a {\em distant parallelism} defined by the frame choice as one of the possible ``metrical connections''.\footnote{So far, the specification of the translational part of the Cartan connection could be left open. It entered the investigation only when  Cartan  characterized the condition of vanishing torsion (below PA).  }
Note that this was 
{\em  no faithful representation} of Weyl's  description of ``metrical'' infinitesimal geomery with its   {\em   specific difference}  between infinitesimal {\em similarities} $\mathfrak g ^{\ast}$(for metrical displacements) and {\em congruences} $\mathfrak g$ (for ``rotations'' at a point).
 In his mode of rephrasing Weyl's idea, Cartan suppressed this specific difference.\footnote{In modern parlance, Weyl's construction is not  described by a principal fibre bundle with structure group $G$, but by a bundle 
 $E \rightarrow M$ with fibre $G$, and the tangent bundle $TM$ associated to it. The structure group of $E$ is $G^{\ast}$ which  operates on $G$ by conjugation. }
 
 In the next step, Cartan considered, like Weyl, the whole collection of ``metrical connections'' derivable from the distant parallel displacement by adding a differential form with values in $\mathfrak g$ to the one given by the frame choice. He then translated Weyl's stipulation of symmetry for an affine connection into his own terms:\footnote{``Par suite ,l'axiome I de M. H. Weyl peut s'ènoncer de la mani\`ere suivante: \ldots'' \cite[173]{Cartan:1923PoS}}
\begin{itemize}
	\item[PA] Whatever the choice of the translational part, among all Cartan connections (values in  $\R^n + \mathfrak g$) with rotational part  formed in this way there is at least one with vanishing torsion (ibid. 173). 
\end{itemize}
As ``second axiome de M. H. Weyl'' he stated the  uniqueness postulate:
\begin{itemize}
	\item[PB] For a metrical connection in the sense above, satisfying postulat PA, the torsion free connection is uniquely determined.  
\end{itemize}

Of course, Cartan shifted again the content of Weyl's principles: He skipped completely the content of the  ``postulate of freedom'' P1  and separated the ``postulate of coherence'' P2 into two parts, existence (PA) and uniqueness (PB) of a compatible torsion free connection. From a mathematical point of view, both modifications made sense (but Cartan wrote as if he did not notice the shifts). Although Weyl considered P1 as philosophically important and although it led to an easy accessible counting argument for the dimension of $G$, it turned out, in the end, as mathematically redundant. Not only did it become superfluous in Cartan's modified space problem;  E. Scheibe later showed that also for the  Weyl's original PoS it was unnecessary \cite{Scheibe:Diss,Scheibe:Spacetime}.\footnote{The dimension condition follows alredy from the existence postulate of a compatible affine connection.} 

The second modification (separation of P2 in two parts, PA and PB) was helpful for reformulating  the existence condition of vanishing torsion in Cartan's own terms  and for structuring the proof of the main theorem of the PoS in the framework of his theory of infinitesimal Lie groups. Already the existence assumption of a compatible torsion free connection (PA) allowed Cartan to select a certain number of simple groups as  possible candidates for $G$ from his general list \cite[181ff.]{Cartan:1923PoS}.\footnote{The list for the homogeneous part $G$ was: $SL(n,\R)$, subgroup of $SL(n,\R)$ keeping a line fixed, $SO(p,q;\R), Sympl(2n,\R)$. }
 Adding  the uniqueness demand (PB) then reduced the allowable groups to the special orthogonal ones:
\begin{quote} {\em The linear group of a non-degenerate quadratic real form ist therefore the only one which satisfies the conditions of the axioms I and II of Mr. H.Weyl} [here abbreviated as PA and PB, ES]. In this way we have  proven, on the basis of these axioms,  that the Pythagorean form of the metric in the universe is necessary.$^{(vii)}$ 
 \cite[192, emphasis in original]{Cartan:1923PoS}
\end{quote}

Geometrically,  Cartan's adaptation of the modern PoS led to the 
\begin{result}
 Infinitesimal congruences in the sense of Cartan's PoS, satisfying postulates PA and PB allow to reconstruct the structure of a Riemannian metric, comparable to the reconstruction of the Riemannian metric  by Helmholtz, Lie, and Killing in the classical PoS.
\end{result}

Cartan seemed to be satisfied,  philosophically, to have found a  characterization of Riemannian geometry among his much more general class of infinitesimal geometries, just like Weyl had been  about his slightly different argument for his gauge generalization of Riemannian manifolds.

Of course this excursion into the PoS, triggered by Weyl's example, contained only a selected aspect of Cartan's own thoughts about the question of how to modify the concept of space in the light of relativity. He pursued a much vaster program for revising the ``classical''  concept(s) of space, oriented at the goal of achieving for  Klein's view of geometries  a similar infinitesimalization (present physicists' language ``localization'') which Riemann had proposed for Euclidean geometry. In his talk at the Toronto ICM Cartan explained this view in general terms
\begin{quote}
\ldots while a Riemannian space does not possess absolute homogeneity, it  possesses a kind of infinitesimal homogeneity; in the immediate neighbourhood of a point it can be assimilated to a  Euclidean space. 
 \cite[85]{Cartan:1924Toronto}$^{(viii)}$
\end{quote}

For the investigation of such infinitesimal laws of homogeneity, Cartan extended the method of {\em moving frames} (``r\'ep\`eres mobiles'') and {\em Pfaffian (differential) forms} taken over from his teacher Darboux, but generalized both in a highly sophisticated way.
Infinitesimal neighbourhoods (modernized $T_pM$)  of a Riemannian manifold  could be  ``assimilated'' (Cartan  in Toronto)  to an Euclidean space, which again could be constructed as a homogeneous space $\E^n \cong L/G$, for example  from the Euclidean isometries $L\cong \R^n \rtimes SO(n,\R)$ and the rotations $G=SO(n,\R)$.

 Taking a consequently  infinitesimal perspective, Cartan preferred  to pass over to the corresponding infinitesimal Lie groups, in later terminology the Lie algebras $\mathfrak l = Lie \, L, \mathfrak g = Lie \, G$; then the infinitesimal neighbourhood could be ``assimilated to''   a homogeneous space of infinitesimal Lie groups\footnote{Cf. fn. (\ref{Lie algebra}).}  
\beq T_pM \cong \mathfrak l / \mathfrak g  \cong \mathfrak{h}\, . \label{Cartan gauge}
\eeq
where   $\mathfrak h$ was here the trivial subalgebra $\R^n$, invariant under the adjoint operation of $\mathfrak g$ (today {\em reductivity}), and $\mathfrak l = \mathfrak h \oplus \mathfrak g$.
% $\mathfrak l \cong \R^n \otimes \mathfrak s \mathfrak o (m,\R)$. 
Similarly for other pairs of groups $L, G$ for which this  condition ((\ref{Cartan gauge}) and reductivity) was  satisfied. Among them Cartan studied for  $L$ the affine, the conformal and the projective groups, while  $G$ was the respective isotropy group of a point. 
Cartan described such an assimilation by introducing increasingly  general (and tricky) {\em (Cartan) frames} (``r\'ep\`eres''). 
For completing the idea of ``homog\'en\'eit\'e infinit\'esimale'' Cartan explained, how the group operation of $L$, respectively its infinitesimal version $\mathfrak l$ binds the infinitesimal neighbourhoods together. To do this he introduced a collection of ``Pfaffian'' (i.e., differential) forms 
$ \omega^i$ ($1\leq i\leq n$) and $ \omega^j_k $ ($ 1 \leq j, k \leq dim \, \mathfrak g$).\footnote{Different from Cartan's notation,  an upper and lower index notation is used here in analogy to the usage in co- contravariant tensor calculus.}
 They parametrized the infinitesimal operations of  $\mathfrak l / \mathfrak g \cong \mathfrak h$ and $\mathfrak g$ with regard to the chosen Cartan frame. 
 They  prescribed  infinitesimal group operations in $\mathfrak l/ \mathfrak g  \cong \R^n$, respectively $\mathfrak g$ 
 for any infinitesimal variation  $\delta x = (\delta x_1, \ldots, \delta x_n)$  of  coordinates of a given point. In slightly modernized terminology
 \beqa (\omega^i) &\quad & \mbox{is a differential 1-form with values in} \; \mathfrak l/ \mathfrak g \cong \R^n \nonumber \\
 (\omega^j_k) &\quad & \mbox{is a differential 1-form with values in} \; \mathfrak g
\eeqa
Together they form what  later would be called a {\em Cartan connection}. The $(\omega^i)$  represent the generalized translational%\footnote{Today sometimes called  contribution{\em transvections}.} 
and the $(\omega^j_k)$ the isotropic (generalized rotational) contribution of the connection.\footnote{Compare the more detailed presentation in Nabonnand, this volume.}

If one follows (integrates) such motions along a path the result is, in general, path dependent. Infinitesimal closed contours lead, like in the case of the Levi-Civita connection of Riemannian geometry, to curvatures $(\Omega^i)$ (translational part) and $(\Omega^j_k)$ (isotropic part).  Cartan was able to express these curvatures by differential 2-forms with values in the respective infinitesimal groups and   depending on the ``sides''  of an infinitesimal parallelogram 
\beqa\Omega ^i   &=& d\omega ^i -  \omega ^k \wedge  \omega _k^{\;i} \, ,      \, \label{torsion} \\ 
    \Omega _i^j   &=&  d\omega _i^{\; j} -  \omega _i^{\; k} \wedge \omega _k^{\; j}     \; .\quad \quad  \ \label{curvature}
\eeqa 
For the translational part of the curvature, $\Omega ^i $, Cartan coined  the term {\em torsion},  the spaces he called {\em non-holonomous}.

During the 1920s  Cartan  studied non-holonomous spaces of  increasingly complex type: 
 \begin{itemize}
\item Poincar\'e group in papers on the geometrical foundation of  general relativity \cite{Cartan:1922[61],Cartan:1923[66],Cartan:1924[69]}. For torsion  $\Omega ^i =0$ such a Cartan space reduced to a Lorentz manifold and could be used for treating Einstein's theory in Cartan geometric  terms,
\item  inhomogeneous similarity group (for torsion $ =0$, this case reduced to Weylian manifolds),
\item  conformal group \cite{Cartan:1922[60]},
\item  projective group \cite{Cartan:1924[70]}.
\end{itemize}

 In this way, Cartan developed an impressive  conceptual frame for studying different types of differential geometries, Riemannian, Lorentzian, Weylian, affine, conformal, projective. All were  enriched  by the possibility to allow for the new phenomenon of torsion, and all  arose from Cartan's unified method of adapting the Kleinian viewpoint to infinitesimal geometry. 
 
 For him, the look at PoS in the sense of Weyl was  only a small subject in a much wieder field. Nevertheless Cartan mentioned it prominently at the end of his Toronto talk:
 \begin{quote} Before I come to the end, I have to indicate the remarkable researches in which Mr. H. Weyl has taken up the old philosophical problem of the problem of space, which was treated formerly  by Helmholtz and Lie; he has it  now  adapted  to the new viewpoints introduced by relativity theory. The concept of group is here again at the base even of the formulation of the problem posed by Mr. H. Weyl. I cannot dream, however, to enter into even a comprehensive discusssion of this important question which would demand a special presentation of its own. \cite[93f.]{Cartan:1924Toronto}
 $^{(ix)}$
 \end{quote}

This could be read as an invitation to 
 mathematicians for taking  the   modern PoS  seriously as a research problem.  Although Cartan himself had  subsumed  it under his more general perspective of non-holonomous spaces and (Cartan-) connections, he did not consider the  PoS  as closed by his short articles at the beginning of the 1920s.

 \subsubsection*{6. Outlook:  mathematical and physical aspects of the   PoS in the second half of the 20th century}
Most {\em mathematicians} of the next (post-Cartan/Weyl) generation learned  the modernized PoS from Cartan, and thus separated from Weyl's gauge geometric perspective (Tits, Freudenthal,  Nomizu, Klingenberg, e.a.). With the rise of the fibre bundle language, even Cartan's geometry tended to be reduced to the homogeneous (Riemannian) standpoint, and its generalization of $G$-structures (on the tangent bundle of $M$). Only at G\"ottingen two young mathematicians (the second one later turned towrd philosophy of science) did research on the specific Weylian versions (first and second) of the PoS keeping close to Weyl's specific geometric contextualization of his first (1919) and second (1921--23) discussions of the space problem \cite{Laugwitz:ARP,Scheibe:Diss}. 

The other authors of the next generation approached the problem from views closer to main developments in the field. J. Tits and H. Freudenthal analyzed  general topological characterizations of the group operation.  They thus continued and expanded Hilbert's characterization of the foundations of geometry by  topological group principles in the notes of 1902/1903 \cite{Tits:PoSI,Tits:PoSII,Freudenthal:ARP1956}. Moreover they 
simplified the group theoretic part of the investigation and updated it from the point of view of  the structure classification of Lie algebras, which became famous in their generation. The topological group aspect of the PoS had been investigated in the decades before by a series of researchers, among them Brouwer (1909/10), R.L. Moore (1919), Cairns (1923), S\"uss (1925-27), Lubben (1928), Kolmogorov (1930), Montgomery/Zippin (1940), K\'er\'ekart\`o (1950).\footnote{List according to \cite[16]{Freudenthal:ARP1960}; see there for bibliographical references.}  
In their topological characterization, Tits and Freudenthal built along the line proposed by   Kolmogorov. For the underlying space $S$ they  assumed no longer the manifold property, but only metrizability, local compactness and local connectedness. The transformation group $G$ had to satisfy some kind of ``rigidity'', completeness, and a separation property for orbits. Already these quite general conditions led to specifications for the group that could be evaluated in the light of the classification theory of simple Lie algebras.\footnote{\cite{Tits:PoSII,Freudenthal:ARP1960}. Freudenthal commented that Hilbert's 3-closedness condition could be considered as a topological kind of ``rigidity''.} 

The differential geometers of the next generation were strongly influenced by the Cartan tradition. They developed a modernized language for infinitesimal connections,  \cite{Ehresmann:1951,Chern:1951,Nomizu:1956} which fitted well to the new theory of fibre bundles established and studied by topologically oriented geometers  \cite{Steenrod:Fibrebundles,Hirzebruch:Methoden}.  Cartan's analysis of the PoS was translated into this newly established setting. 
In this framework  Cartan's PoS was generalized to any type of $G$-bundle of frames $E_G \rightarrow M$ in the tangent bundle over a (differentiable) manifold $M$, for any closed subgroup $G \subset GL(n,\R)$. Klingenberg then posed the generalized PoS in the form:
\begin{quote}
We now ask for which closed subgroups $G$ of $GL(n,\R)$ there always exists a canonical linear connection in $E_G$, i.e., a connection in some  way specified by the structure of $G$ \ldots \cite[300]{Klingenberg:ARP1956}$^{\mbox{(x)}}$
\end{quote}
He  showed that, e.g., for $2n$-dimensional almost complex manifolds $M$ a result analogous  to Cartan's holds: If for a  closed subgroup $G\subset GL(n,\C)$ any $G$-frame bundle $E_G$ over $M$ allows exactly one connection for which the mixed terms of the torsion vanish, then $G$ is a real form of $GL(n,\C)$ and vice versa.\footnote{Klingenberg generalized a  result  by Chern:  There is exactly one connection with vanishing mixed torsion terms  for the unitary group $U(n)$ and an almost complex manifold \cite{Chern:1946}. In this way Weyl's demand of torsion free connection was weakened to the vanishing mixed torsion terms.}
 Thus Cartan's PoS developed into a wider research field later called the study of {\em $G$-structures} on $M$ \cite[288ff.]{Kobayashi/Nomizu}.

 Parallel to this mathematical research on the PoS, following the lead of Cartan and Weyl, {\em physicists} developed their own understanding of infinitesimal, local and global structures of spacetime and the dynamical fields upon it. Here we can, of course, only scratch upon the surface of this complex and manysided history, which covers large parts of fundamental physics of the second half in the 20th century. We can only roughly indicate  two developments in which physicists started to see the role of ``localized'' symmetries  in spacetime (as they used to call  point dependent infinitesimal operations of Lie groups)  or acting in  dynamical spaces above it ({\em internal} spaces). Both are most closely  related to Weyl's and Cartan's treatment of the PoS. 
 
 These developments come  from the two branches of physical theory which, about the middle of the last century,   started to be built about gauge structures and  turned the latter into  a central paradigm in  fundamental physics of the outgoing 20th century.  The most spectacular development in this direction arose 
  from the growth of gauge field theory in elementary particle physics. It was  opened up by the work of \citeasnoun{Yang/Mills} and others, less well known but earlier and more general, Utiyama who sketched an integrated gauge panorama for both, gravity and nuclear fields.\footnote{See \cite{Pickering:Quarks,RiseSM,%
 ORaifeartaigh:Dawning,ORaif/Straumann}.}
 The other one consisted  in the deepening usage of Cartan geometric structures  in  gravity theory.\footnote{See \cite{Blagojevic/Hehl}.} 
 
 At the end of  section 3, we saw a widening gap in general relativity between mathematical infinitesimal symmetries and their physical interpretation, as far as energy conservation was concerned. That stood in clear contrast to special relativity, where the generalized ``motions'' (automorphisms of Minkowski space) not only could be used to describe transformations of rigid measuring devices between different inertial systems. On a much deeper level, these automorphisms   turned out  to be the structural origin of  conservation laws, in particular of energy-momentum. In GR the question of how to treat energy conservation was more involved and remained a challenge for a long time. In general relativity, infinitesimal translations could be modelled in two ways: by infinitesimal diffeomorphisms of the manifold or  through infinitesimal translations  in the infinitesimal homogeneous spaces of  Cartan geometries. Both ways, the variation of a  dynamical Lagrangian led to  
 energy momentum currents,  the  {\em canonical energy momentum}  $\mathfrak t_{ij}$ for variation with regard to 
 infinitesimal translations, while  the {\em dynamical energy momentum}, $T_{ij}$,   was  derived from varying  the gravitational potential (in Einstein gravity  $g_{ij}$, in the Cartan geometric generalizations the respective connections).   At first it was a problem how to make these two expressions consistent.  First answers were given in the 1940s by Belinfante, Rosenfeld (see end of section 3).

  D. Sciama and T. Kibble, both working in London,  at King's College and   Imperial College respectively, investigated different but methodologically closely  related approaches for a consequent ``covariantization'' of a special relativistic field theory, i.e., a generalization and import into general relativity. Their idea was to ``localize''  the symmetries of Minkowski space, i.e., the Poincar\'e group \cite{Sciama:1962,Kibble:1961}. In this enterprise they  made use of  the basic structures of Cartan geometry, although   with quantities expressed in classical tensor calculus, without explicit recourse to Cartan's calculus of differential forms.\footnote{Both quoted \cite{Weyl:1929Dirac} which used certain aspects of Cartan's work and anticipated a gauge theoretic formulation of general relativity. Sciama referred to Weyl as his source for the use of  ``vierbeins''. He did not mention   Cartan  nor his differential forms. Notice the paradoxical historical analogy with the  transmission of Weyl's gauge theory of 1929 to Yang and Mills: Yang and Mills did not know of Weyl's work but absorbed its  basic idea in a simplified form through Pauli's article on wave mechanics and generalized it.}
  
 It did not take long until Cartan's geometric  methods, sometimes combined  with tensor calculus, were introduced  into this research program. A. Trautman and F. Hehl, with their Warszaw respectively Cologne collaborators, seem to have been among the earliest to do so \cite{Trautman:1973,Hehl_ea:1976}.\footnote{More detailed historical research on this question is necessary. At the moment the best information on this complex can be found in the  reader  \cite{Blagojevic/Hehl} which contains extremely helpful commentaries. For ex-post surveys see  \cite{Trautman:2006,Hehl:Dennis}. } 
Their work  showed that Cartan geometry offered a tailor-made geometric framework for  infinitesimalizing  (``localizing'' in the language of physicists) the currents known from Minkowski space and special relativity.\footnote{That had not been done in the earlier generation, maybe because
Cartan himself, as well as Einstein, Weitzenb\"ock and other physicists of the first generation of relativists, used Cartan geometry  to rewrite Einstein gravity in terms of  {\em  distant parallelism}, see \cite{Sauer:Fernpar,Goenner:UFT}.
 In this setting the deviation of flat space was encoded in the torsion part of curvature, while rotational curvature was set to zero.}

In the Sciama-Kibble approach, Cartan geometry was allowed to have both, rotational curvature and torsion, and was additionally endowed   with a Riemannian metric ({\em Riemann-Cartan} geometry).
 The gravitational Lagrangian was  modelled as closely as possible to the Hilbert action of Einstein gravity, but the geometrical structure was fully Cartanian (with Poincar\'e group). As the infinitesimal neighbourhoods of the  (spacetime) manifold were identified with  homogeneous spaces constructed from the Poincar\'e group, i.e., affine rather than vector spaces,
 point dependent infinitesimal translations could be handled in  analogy to gauge transformations in elementary particle physics.\footnote{Transformations of these geometric  {\em Cartan gauges}, in the terminology of \cite{Sharpe:DiffGeo}, consist    of two components, one translational transforming as a differential form,  the other rotational transforming as one expects from gauge connections. }
  In the dynamical equations of this approach (a modified Einstein equation and another one regulating  spin-torsion coupling)   the canonical energy-momentum appears on the right hand side of the modified Einstein equation (similarly the torsion equation has the spin current as source). It  is  conserved, due to the underlying Cartan geometric symmetry structure \cite{Sciama:1962,Hehl_ea:1976,Trautman:2006,Hehl:Dennis}.
 This was a striking {\em structural insight}, independent of the  empirical relevance of the spin-torsion extension.
  Moreover, but far beyond  the scope of this contribution, the program was extended to studying more general Lagrangians (e.g., quadratic in curvature and torsion) and more general Cartan geometries (in particular { metric affine} Cartan geometry).\footnote{For these developments see \cite{Blagojevic/Hehl}.}

  In this way,   the central  guiding idea of Cartan in the development of his  ``non-holonomic'' (Cartan) spaces was realized not only in geometrical but also in dynamical terms. In his Toronto talk Cartan had put his view like this:
 \begin{quote} \ldots it is nothing but the development of relativity theory,  bound by the paradoxical obligation to interpret, in a non-homogeneous universe, the numerous experiences made by observers which believe in the homogeneity of this universe, which allows to  fill up, at least partially, the trench which separates Riemannian spaces from Euclidean space. \cite[86]{Cartan:1924Toronto}$^{(xi)}$
 \end{quote} 
 The Sciama-Kibble-Hehl-Trautman e.a. research program underpinned  the   dynamical nature of the infinitesimal translations component of what Weyl had called the physical automorphisms of general relativity (see end of section 4) and gave it a  clearer  mathematical expression in the ``localized'' Poincar\'e group.
 
 Also in the 1970s,  Weyl's original gauge geometry with infinitesimalized scale invariance was taken up again by Omote, Utiyama e.a.in Japan    and independently by Dirac, Ehlers, Pirani and Schild and others in Europe.\footnote{\cite{Omote:1971,EPS,Dirac:1973,Omote:1974,Utiyama:1975I,%
Utiyama:1975II}.}
 This retake was no longer motivated by Weyl's original hope for a geometrically unified field theory of electromagnetism and gravity, it rather arose from an attempt to establish a bridge between a scalar field extended gravity  theory of Jordan-Brans-Dicke type and cosmology (Dirac), or to the nuclear forces the theory of which still lay in the dark (Utiyama).\footnote{Dirac, surprisingly, stuck to the obsolete interpretation of Weyl's scale connection as the potential of the electromagnetic field, the other authors followed different lines. For  a recent critical acclaim of the foundational oriented paper of Ehlers/Pirani/Schild see   \cite{Trautman:EPS}.}   
 
 Weyl's scale gauge geometry  offered a chance to build a conceptual bridge between gravity and the conformal structures of fundamental forces, in  as much  as they were considered under abstraction from mass (with ``massless'' fermionic and bosonic fields). This chance did not materialize in the short run. One reason, perhaps the crucial one, was  the unclear role of the scale covariant scalar field of Weyl-Omote gravity, which remained a riddle and open for speculations in different directions. So the approach did not lead to  results  as definite  as the torsion-spin coupling and the clarification of the source currents in Einstein-Cartan gravity. In fact, scalar fields remained  hypothetical until recently (2012) also in elementary particle physics,\footnote{\cite{ATLAS:Higgs2012,CMS:Higgs2012}} although a specific one, the {\em Higgs field}, had become a central gadget already with the rise of the new fundamental theory of matter which developed in the 1970/80s.\footnote{For  attempts to make use of Weyl geometry in relation to fundamental forces see among others \cite{HungCheng:1988,Drechsler:Higgs,Nishino/Rajpoot:2004} or, more coceptually, \cite{Scholz:Dennis}. } 
 
 In the meantime, gauge ideas also gained ground in electromagnetism and elementary particle physics. ``Localized'' group extensions of the Lorentz group were explored, in the hope to find mathematical symbolizations of the strong, later also the weak,  interaction, or even a combined ({perhaps even unified}) electro-weak dynamics.\footnote{\cite{ORaifeartaigh:Dawning,ORaif/Straumann,Karaca:Higgs}.}
  In these approaches the gauge idea was stripped from the explicit relationship to gravity, even though by its very nature gauge fields shared ``the character of general relativity'',  as Weyl had expressed it. Such a link  was not obvious, but it was around and served as a  heuristic idea to many of the actors of the post-Weyl/Cartan generation. 
   Originally, these attempts were speculative, tentative and loaden with difficulties. Before the  mid-seventies only  few physicists  would have expected the striking success which this theory form acquired with the rise of the standard model of elementary particles.
   The main task was to bring order into the plethora of new particles or resonances, to find principles to structure them, and to look for overarching theoretical explanations of the dynamics. Groups, weight diagrams of group representations, and infinitesimalized operations (gauge structures) played an increasing role in the 1960s and contributed piece by piece to forming  the  model class which agglomerated and became known, from 1974 onward, as the  standard model of electroweak ($SU(2)\times U(1)$) and chromodynamic/strong ($SU(3)$) interactions. Leptons and quarks became now  the fundamental fermionic constituents of matter.\footnote{\cite{Pickering:Quarks,RiseSM}.}
   
   A crucial contribution to the final success of the standard model  was the realization of the renormalizability of gauge theories. That was  an important precondition for drawing quantitative consequences of the theories by perturbative calculations of the electroweak interaction or, for QCD, by   lattice approximations.\footnote{\cite[341f.]{Kragh:Generations}} For a while, renormalizability was considered as a most important selection criterion for theory approaches. 
  Crucial relations for renormalizability (Ward-Takahashi relation for electromagnetism and the Slavnov-Taylor identities for the non-abelian generalizations) turned out, in due time, as quantum field consequences  of infinitesimal gauge symmetries \cite[chaps. 11.2, 12.4]{Itzykson/Zuber}, although this was not at all clear to the inventors of the first renormalizability proofs \cite{tHooft:HistRenormalizability}.\footnote{I thank A. Borrelli for pointing this out to me. For a historical discussion see \cite[sect. 6]{Zuber:symmetry}.}
  
  In the last third of the 20th century ``localized'' symmetries started to play a crucial role  in both branches of the relativistic field theories,  partly with the intention to trace conserved quantities (gravity theory), partly because of their structural importance for renormalizability (elementary particle theory).   In both cases the infinitesimal symmetry structures conceived by Weyl and Cartan in the 1920s  spread from the analysis of spacetime to the theory of matter and started to give them a common conceptual framework.  In this mutual assimilation, features  carrying ``the character of general relativity'', to put it in Weyl's terms (end of section 4), have been imported into the basic structures of elementary particle physics; on the other hand the gauge structures of gravity have become manifest. Weyl's modified $U(1)$ gauge structure, mimicking characteristic features of infinitesimal symmetries explored in his analysis of the PoS was an inspiring background for this development (indirect for Yang and Mills, direct for Sciama and Kibble), while Cartan's generalized geometries delivered the appropriate geometric framework for a broader ``gauging of gravity''.\footnote{\cite{Blagojevic/Hehl}, here in particular p. 179.}  Seen from this  vantage point,    the fundamental theories of the structure of matter and of gravity are  already permeated by a {\em coherent conceptual approach}, even though a  { unification} of gravity  physics and relativistic quantum field theory in the strong sense, hoped for by many physicists and some philosophers, is still lacking and may  continue to be so.\footnote{The question, whether a quantization of the geometrical degrees of freedom is advisable, is still wide open. A closer inspection of inter-theory relations is, of course, needed and under way. For a lucid survey of the present situation in this regard, see \cite{Fredenhagen_ea:Where_we_are}.}   
   \\[5em]

%\newpage
\subsubsection*{Endnotes (citations in orginal language)}
\begin{enumerate}
\item[\quad (i)] Der zehngliedrigen Gruppe der ``Bewegungen'' im $(x, y, z, t)$-
Raum mit der Ma\ss{}bestimmung:
\[  ds^2 = dx^2 +dy^2 +dz^2-dt^2 \]
entsprechend, gelten f\"ur den ganzen K\"orper 10 dem Schwerpunkt-,
Fl\"achen- und Energieprinzip der gew\"ohnlichen Mechanik
analoge S\"atze. \cite[511]{Herglotz:1911}

	\item[\quad (ii)] Das Beispiel des Raumes ist zugleich sehr lehrreich f\"ur diejenige Frage der Ph\"anomenologie, die mir die eigentlich entscheidende zu sein scheint: inwieweit die Abgrenzung der dem Bewu\ss{}tsein aufgehenden Wesenheiten eine dem Reich des Gegebenen selbst eigent\"umliche Struktur zum Ausdruck bringt und inwieweit an ihr blo\ss{}e Konvention beteiligt ist. \cite[133f.]{Weyl:RZM4}
	\item[(iii)] An die Stelle der von Helmholtz geforderten Homogenit\"at des metrischen Feldes ist die M\"oglichkeit getreten, {\em im Rahmen der feststehenden Natur der Metrik das metrische Feld beliebigen virtuellen Ver\"anderungen zu unterwerfen}. 
	 \cite[46, Hervorh. im Original]{Weyl:ARP1923}
	 \item[(iv)] Ebenso zwingend wie die alte Auffassung einer a priori dem Raume innewohnenden und von der Materie unabh\"angigen metrischen Struktur zur Helmholtz-Lieschen Kennzeichnung der Drehungsgruppe f\"uhrt, l\"a\ss{}t die moderne, von der Einsteinschen Relativit\"atstheorie ausgebildeten Auffassung, nach welcher die  Ma\ss{}struktur ver\"anderungsf\"ahig ist, und in kausaler Abh\"angigkeit von der Materie steht, eine andere Eigenschaft der Drehungsgruppe als die entscheidende erkennen \ldots \cite[348]{Weyl:1923Drehungsgruppe}
\item[(v)] Unsere Partie auf dem Schachbrett des Matrizenschemas ist zu Ende gespielt. (ibid., 372)
\item[(vi)] Cela pos\'e, la signification du premier axiome de M. H. Weyl est la suivante. Choisissons arbitrairement en chaque point $P$ de l'espace l`{\em orientation} du groupe $G$ des rotations, c'est \`a dire les coefficients $a_{ij}$  \ldots [which describe the operation of the group, ES]. Un tel choix \'etant fait, la connexion m\'etrique de l'espace est d\'etermin\'ee: l'une des correspondances par congruences    \ldots est celle pour laquelle les deux vecteurs correspondants ont {\em m\`emes} composantes \ldots \cite[172]{Cartan:1923PoS}
\item[(vii)] {\em Il n'y a donc que le groupe lin\'eaire d'une forme quadratique r\'eelle non d\'egener\'ee qui satisfasse aux conditions pos\'ees par les axiom I et II de M. H. Weyl} [here abbreviated as PA and PB, ES]. Et c'est ainsi qu'est d\'emontr\'ee, en partant de ces axiomes, la n\'ecessit\'e de la forme pythagorienne de la m\'etrique d'Univers. \cite[192, emphasis in original]{Cartan:1923PoS}
\item[(viii)] \ldots si un espace de Riemann ne poss\`ede pas une homog\'en\'eit\'e absolue, il poss\`ede cependant une sorte d'homog\'en\'eit\'e infinit\'esimale; au voisinage immediat d'un point donn\'e il est assimilable \`a une espace euclidien. \cite[85]{Cartan:1924Toronto}
\item[(ix)]  Je ne puis enfin terminer sans signaler les remarquables recherches dans lesquelles M.H. Weyl a repris l'ancien probl\`eme philosophique de l'espace, trait\'e autrefois par Helmholtz et Lie, pour adapter aux points de vue nouveaux introduits par la th\'eorie de relativit\'e; la notion de groupe est, l\`a encore, \`a la base de l'\'enonc\'e m\^eme du probl\`eme pos\'e par M.H. Weyl. Mais je ne puit songer \`a entrer dans l'exposition, m\^eme sommaire, de cette importante question, qui exigerait \`a elle seule une conf\'erence speciale. \cite[93f.]{Cartan:1924Toronto}
\item[(x)] Wir fragen nun, f\"ur welche abgeschlossenen Untergruppen $G$ von\\ $GL(n,\R)$ es in einem $G$-B\"undel $E_G$ stets einen kanonischen, d.h. durch die Struktur von $E_G$ irgendwie ausgezeichneten,  linearen Zusammenhang gibt; \ldots 
 \cite[300]{Klingenberg:ARP1956}.
 \item[(xi)] Or, c'est le d\'eveloppement
m\^eme de la th\' orie de la relativit\'e, lié par l'obligation paradoxale d'interpr\'eter
dans et par un Univers non homog\`ene les r\'esultats de nombreuses exp\'eriences
faites par des observateurs qui croyaient \`a l'homog\'en\'eit\'e de cet Univers, qui
permit de combler en partie le foss\'e qui s\'eparait les espaces de Riemann de
l'espace euclidien. \cite[86]{Cartan:1924Toronto} \\[3em]

\end{enumerate}

\subsubsection*{Glossary}
The following explanations do not provide technically precise definitions. They only serve as a  short heuristic,  explanatory guide to some of the technical terms used in the article. For more details  see {\em Wikipedia} or,  even better, the corresponding literature. 
\begin{itemize}

	\item {\em Cartan space} --  a space (differentiable manifold) $M$ which ``looks infinitesimally'' like a Klein space, i.e. a homogeneous space $G/H$ of a ($\rightarrow$) Lie group $G$ factored by a closed subgroup $H$. (Paradigmatic example: $G$ the real affine group in $n	$ dimensions, $H \cong Gl(n,\R)$ the homogeneous linear transformations, $G/H$ the cosets of the homogeneous transformations (the isotropy group of $O$). To each such coset corresponds a translation and vice versa). ``Infinitesimally'' means that the tangent space $T_pM$ at any $p$  can be identified with  the quotient of the corresponding Lie algebras $\mathfrak{g}/\mathfrak{h}$. Moreever, every point $p$ is linked to an infinitesimal close one $p'$ by a ($\rightarrow$) (Cartan-) {\em connection} which tells how to compare, or to transfer elements from  $T_pM \sim \mathfrak{g}/\mathfrak{h}$  between $p$ and $p'$.
	
		\item {\em conformal structure} --  encodes the angle information of ($\rightarrow$) Riemannian geometry.  In Lorentzian geometry it is  characterized by the light cones at every point and thus, physically speaking, the causal relationships in spacetime ({\em causal structure}). Formally it is  characterized by a class  $[g_{ij}]$ of equivalent metrics, $g_{ij} \sim g'_{ij}$, with $g'_{ij}(x) = \lambda(x) g_{ij}(x)$ (at all points $x$ of a manifold) for some strictly positive function $\lambda:M\rightarrow \R^+$. 
 
	\item {\em connection} (in the sense of differential geometry) -- a symbolic gadget telling how, in a differentiable manifold $M$, elements in spaces $F_p$ adjoined to  points $p\in M$ (called the fibres over $p$, ($\rightarrow$ fibre bundle)) can be compared, or transferred, from a point $p$ to an infinitesimal close one $p'$. More technically, the ``gadget'' can be characterized by a differential 2-form with values in a Lie algebra $\mathfrak{g}$ of a group $G$ operating typically in every $F_p$. 
	In  case of the tangent structure  $F_p=T_pM$ of a Riemannian manifold $M$,  the connection can be given by a 2-form with values in $\mathfrak{g}=gl(n,\R)$ (i.e. $n\times n$ matrices) or, in a reduced form,  in $so(n,\R)$. Best known is the classical {\em Levi-Civita connection} of a Riemannian manifold, given by  the Christoffel symbols $\Gamma_{ij}^i$  which can be expressed by the metrical coefficients $g_{ij}$ and its partial derivatives. In a ($\rightarrow$) Cartan space modelled on a Klein space of type $G/H$ the connection has values in the Liealgebra $\mathfrak{g}$ of the ``large'' group $G$. Usually (in the ``reductive'' case) $\mathfrak{g}$ can be decomposed into a component with values in $\mathfrak{h}$, the Liealgebra of the ``small'' group $H$, and a  component with values in a complementary subalgebra $\mathfrak{l}$, where $\mathfrak{g}\cong \mathfrak{l} \oplus \mathfrak{h}$.
	
	\item {\em covariant differentiation} -- a differentiation procedure defined with regard  to a ($\rightarrow$) connection. It operates on functions (fields) with values in the spaces $F_p$ (fibres of some kind) which are ``transferred'' by evaluating  the connection. It is called ``covariant'' because the derivative contains a co-tensor component more than the function/field that one derives. In index notation: if $X$ is the field which is to be covariantly derived (with regard to some connection), then the covariant derivative $DX$ has components $D_iX$ which transform ``covariantly'', i.e. like covectors (vectors of the dualized tangent space) in the indices $i$. A field is considered to be {\em parallel} with regard to a connection, if its covariant derivative vanishes everywhere in its domain of definition. A field is considered to be {\em parallel transported} along a curve $\gamma$, if the restriction of the field to $\gamma$ is parallel. Parallel transport along different paths with the same initial and end points usually lead to different results.

		\item {\em curvature} of a connection (at a point $p$) --  transformation in $g$ induced by  parallel transport with regard to a ($\rightarrow$) connection about  infinitesimally small closed loops starting and ending at $p$.  In case of the Levi-Civita connection, one gets the {\em Riemann curvature}; in case of a Cartan connection the 
		{\em Cartan curvature}. In the latter case,   the curvature can often  be decomposed (e.g. in the case $G=$ affine group, more generally in the ``reductive'' case)  into a (generalized) {\em rotational part} with values in  $\mathfrak{h}$  and a ``translational'' one with values in $\mathfrak{g}/\mathfrak{h}\cong \mathfrak{l}$. Because of Cartan's  physical interpretation  of the last term by means of  generalized elastic media (Cosserat media), this term is called  {\em torsion}.

	\item {\em distant parallelism} -- a gadget which allows to compare directions in tangent spaces $T_pM$ and $T_qM$ for any two points $p, q$ of a differentiable manifold $M$. More technically, a distant parallelism of manifold of dimension $n$ boils down (is equivalent) to specifying $n$ nowhere vanishing vectorfields which are everywhere (at every point) linearly independent. 
	Because of topological obstructions, this is often not possible. If it is, the manifold is called ``parallelizable''. In works of mathematical physics (unified field theories) of the 1920/30s, the existence of  a distant parallelism was usually assumed to be unproblematic (which it is, if one restricts to local considerations).   
	
	\item {\em energy-momentum  (e-m) tensor} -- a tensor, often denoted by $T =(T_{ij})$, which contains all quantitative information on energy density, pressure, energy flow and momentum flow of a physical field. In  general relativity it forms the right hand side of the Einstein equation and can be formulated in different non-equivalent forms. In most cases it is derived from  varying the Lagrangian density of matter terms with regard to the metric ({\em Hilbert energy momentum tensor}). Another one (or better a whole class) arises as a current of ``infinitesimal translations'' considered as  ($\rightarrow$)   Noether symmetries ({\em canonical e-m tensor}). If it is clear which choice is correct for describing a dynamical situation (in most cases a gravitational one), the specified tensor is called the {\em dynamical e-m tensor}. Often the Hilbert tensor is considered as the dynamical one. 
		
	\item {\em fibre bundle} -- a big manifold (the ``total space'') ${P}$ lying over  over a small one (the ``base'') $M$ by a map   $\pi:\, {P}\rightarrow M$. All counterimages of base points $\pi^{-1} (p)$, the {\em fibres},  are assumed to be isomorphic among each other and to a {\em standard fibre} $F$. Moreover a group $G$ operates on all fibres in a way which can  locally  be described in a simple form through the operation of $G$ on $F$ in a ``local trivialization'' which is comparable to a  coordinate choice in a manifold. In the physics context, one often speaks of chosing a ``gauge'' rather than choosing a local trivialization.   Different local trivializations are related by point-dependent operations of the group $G$; changes of trivializations lead to ($\rightarrow$) gauge transformations.
	$G$ is  called the {\em structure group} of the bundle. In the special case $F=G$ one speaks of a {\em principal fibre bundle}. If a ($\rightarrow$)  connection is given in a principal bundle with group $G$, any fibre bundle with structure group $G$, a so-called {\em associate bundle} inherits a connection.	
	
	\item {\em gauge transformation} -- a couple of point dependent transformations in  $G$ operating on a ($\rightarrow$) fibre bundle with structure group $G$,  endowed with a ($\rightarrow$) connection. The functions with values in the fibres (mathematically speaking the ``sections'', physically speaking vector, tensor, or spinor fields) are transformed by appropriate representations of the group; the connection (which can locally be represented by a $\mathfrak{g}$-valued differential 2-form on the base space) by a more tricky type of transformation. These transformations are often called ``gauge transformation of first kind'', respectively of ``second kind''. The gauge transformations are, by definition, point  dependent and operate ``in the infinitesimal''. In the physics literature they are usually called {\em local transformations} (of $G$) in contrast to {\em global} ones, by which $G$ operate  constantly on the fibres, respectively the sections. An {\em active}  gauge transformation  changes the sections/fields and the connection, while a  {\em passive} gauge transformation only changes the representation of the  fields and the connection with regard to a local trivialization, respectively a choice of gauge. Thes two aspects are  comparable to the complementarity of endomorphisms and coordinate transformations for vector spaces, or (local) diffeomorphisms versus coordinate change for manifolds. 
	
	\item {\em Lie group} -- a group $G$ which has  a continous structure given by a differentiable real or complex manifold. Algebraic group  operations have to be differentiable. Infinitesimal operations can be represented on the tangent space  $T_eG$ of the neutral element $e \in G$; they form the {\em Lie algebra} (in older literature ``infinitesimal Lie group''). If a Lie algebra cannot be  be decomposed non-trivially, it is called {\em simple}; simple Lie algebras are the building blocks of the more complicated ones. Simple 	Lie algebras are more easily classified over the complex numbers than over the reals, i.e. as algebras over  $\C$. Thus arises an intriguing interplay of complexification of real Lie algebras and search for real subalgebras of complex ones, which are called {\em real forms} of the complex algebra.

	\item {\em Noether theorems} -- two theorems established by Emmy Noether dealing with consequences of dynamical variation problems with Lagrange densities invariant under constant (global) or point dependent (local) symmetries with regard to continuous groups ($\rightarrow$ Lie group). They are  called Noether theorem I and II,  respectively. If a Lagrangian depending on fields $X$ is invariant under constant symmetries of a Lie group of dimension $n$,  $n$ quantities can be formed from the fields, which remain constant under the time flow of the solutions to the variation problem. They are called {\em Noether charges}. In classical, and in special relativistic, mechanics energy, momentum, and angular momentum are such Noether charges; the corresponding conservation laws are consequences of the Noether thm. I. The second Noether theorem is formulated in a very general form. It deals with  point dependent symmetries and gives partial differential equations ({\em Noether relations}) for each group generator.  Noether thm. II, applied to a Lagrangian density which is invariant under local symmetry operations of a Lie group of dimension $n$, allows to derive, for each Noether relation, a vector (density) field $\mathfrak{j}^k$ which satisfies a divergence relation of type $ \sum_k \partial_k \mathfrak{j}^k = 0$. These are called {\em Noether currents}. 	
	In the context of general relativity  Noether currents with regard to localized (``infinitesimal'') translations have been discussed as candidates for a generalized type of energy and momentum (current) conservation; but problems for a meaningful physical interpretation are considerable (see main text). In gauge field theories of the Yang-Mills type conserved currents are easier to handle and turn out to be important for quantization.  
	
			\item {\em Riemannian metric} -- a metrical concept in differential geometry which generalizes the differential expressions for measurements on curved surfaces embedded in Euclidean 3-space to  $n$-dimensional manifolds $M$. If local coordinates in $M$ are given by a set of real valued functions $x= (x_1, \ldots,  x_n)$ the  ``line element'' $ds$  is given by the formula $ds^2 = \sum_{i,j=1} ^n \, g_{ij}(x_1, \ldots x_n)\, dx_i dx_j$, where $(g_{ij})$ are the entries of a positive definite symmetric, point dependent matrix (two-form). In general relativity notation is usually  simplified by using lower and upper indexes according to the symbolic function of the expressions (vectors,  covectors etc.), and the Einstein summation convention. The coordinate functions are then better denoted by $x^1, \ldots, x^n$ (although they are, of course, no vector components) and the squared line elements become $ds^2 = g_{ij}dx^i dx^j$. A Riemannian metric allows to introduce shortest lines,  {\em geodesics}, and goes in hand with a uniquely determined affine ($\rightarrow$) connection, the Levi-Civita connection, with regard to which parallel displacement, ($\rightarrow$) covariant derivation and (Riemannian) curvature ($\leftarrow$) are  defined.   Because of the role of Minkowski space in special relativity, it becomes  very natural in general relativity to consider  a non-degerate symmetric two-form $g = (g_{ij})$ which is of Lorentz signature, i.e. which after diagonalization and normalization has the form $g= diag(+1, -1,-1,-1)$ or the other way round,  $g= diag(-1, +1,+1,+1)$. Then one speaks of a {\em Lorentzian manifold} or, in more general cases ($g$ non-degenerate, symmetric), of {\em pseudo-Riemannian manifolds}. Unique determination of affine connection, and covariant derivative and curvature remain like in Riemannian geometry.

	\item {\em Weylian metric} (gauge metric) -- a generalization of ($\rightarrow$) Riemannian metrics proposed by H. Weyl in 1918 in order to ``localize'' comparability of length measurement (and of related quantities). Weyl postulated a specification of a ($\rightarrow$) conformal class of (pseudo-) Riemannian metrics $[(g_{ij})]$ (to allow for  length comparison of vectors at the same point) and, in addition, of a class of differential one-forms $\varphi = \sum_i \varphi_i dx^i$ (for length comparison between infinitesimal close points). Any one-form $\varphi$ associated to the choice of a representative $g_{ij}$ of the conformal class of metrics plays the role of a {\em length} or {\em scale connection}. It allows  for ``length transfer'', i.e. comparison between units and length measurements at infinitely close points, and, after integration along paths, for length comparison of vectors at  finitely distant points $p, q \in M$. A choice of the representative can be interpreted as a point dependent {\em gauge} of the units of length. Change of the representative $g_{ij}(x) \rightarrow g'_{ij}(x)= \Omega(x)^2g_{ij}(x) $ (with nowhere vanishing real valued function $\Omega$) corresponded to a change of gauge in the literal sense. Such a change is accompanied by a change of the representative of the length connection $\varphi \rightarrow \varphi' =  \varphi  - d \log \Omega $. This was the first ($\rightarrow$) {\em gauge transformation} explicitly considered  in the mathematical literature. In the terminology of  ($\rightarrow$) fibre bundles it is a very simple connection with values in the trivial Lie algebra of the commutative group $(\R^+, \cdot)$. A Weylian metric has a uniquely determined  compatible ($\rightarrow$) affine connection and covariant derivative. They lead to ($\rightarrow$) curvature concepts similar to those of Riemannian geometry. Moreover there arises a curvature $f= (f_{ij})$  of the length connection. It is easy to calculate:  $f= d \varphi = f_{ij}\, dx^i \wedge dx^j$, $f_{ij} = \partial_i \varphi_j - \partial_j \varphi_i$. If the length curvature vanishes, $f=0$, the Weylian metric can locally be gauged to the form of Riemannian geometry (i.e. $\varphi' = 0$).

\end{itemize}
%\newpage

%\nocite{Weyl:GA}

% \bibliographystyle{apsr}
 % \bibliography{a_lit_hist,a_lit_mathsci}

\end{document}